\documentclass[12pt]{article}
\usepackage{amsthm, amsmath,amssymb,amsfonts}
\usepackage{graphicx}
\usepackage{psfig}
\usepackage{graphics}
\usepackage{amsmath}
\usepackage{amssymb}
\usepackage{epsfig}
\usepackage{epstopdf}
\usepackage{color}

\newtheorem{theorem}{Theorem}[section]
\newtheorem{proposition}[theorem]{Proposition}
\newtheorem{corollary}[theorem]{Corollary}
\newtheorem{lemma}[theorem]{Lemma}
\newtheorem{remark}[theorem]{Remark}
\newtheorem{definition}[theorem]{Definition}
\newtheorem{example}[theorem]{Example}

\def\para{\vspace{2mm}}

\def\gcd{{\rm gcd}}

\def\ord{{\rm ord}}

\def\Re{{\rm Re}}
\def\Im{{\rm Im}}

\begin{document}
\title{Asymptotic Behavior of an Implicit Algebraic Plane Curve}
\author{Angel Blasco and Sonia P\'erez-D\'{\i}az\\
Departamento de F\'{\i}sica y Matem\'aticas \\
        Universidad de Alcal\'a \\
      E-28871 Madrid, Spain  \\
angel.blasco@uah.es, sonia.perez@uah.es
}
\date{}          % Enter your date or \today between curly braces
\maketitle

\begin{abstract}
In this paper, we introduce the notion of infinity branches as well
as approaching curves. We present some properties which allow us to
obtain an  algorithm that compares the behavior of two implicitly defined
algebraic plane curves at  the infinity. As an important result, we prove that if two plane algebraic curves have the same
asymptotic behavior,   the Hausdorff distance between them is
finite.
\end{abstract}

{\bf Keywords:} Implicit Algebraic Plane Curve; Infinity Branches; Convergent Branches; Asymptotic Behavior; Approaching Curves

\section{Introduction}

Unirational algebraic varieties, play an important role in the frame of practical
applications (see
 \cite{HSW} and \cite{HL97}).  In particular, many authors have studied different problems related to plane algebraic curves that are defined implicitly
 (see e.g. \cite{SWP} and \cite{Walker}). In this paper, we deal with the notion of infinity branches which  is a very important tool to analyze the
 behavior of an implicitly defined algebraic plane curve  at    the
infinity.  For instance,  determining the infinity branches of an implicit real algebraic plane curve is an important step in sketching its graph as well as in studying  its topology (see e.g. \cite{Gao}, \cite{lalo}, \cite{Hong} and \cite{Zeng}).

\para

 Intuitively speaking, the infinity branch of a real plane algebraic
curve reflects the status of a curve at the points with sufficiently
large coordinates. An infinity branch is associated to
a projective place centered at an infinity point, and it can be {\it
parametrized}  by means of  Puiseux series. We show how
to obtain this parametrization.

\para

The concept of infinity branch  allows us to introduce the notion of
convergent branches  and approaching curves. Intuitively speaking,
two infinity branches converge if they get closer  as they tend to infinity.
This notion allows us to analyze whether two given implicit algebraic plane  curves
approach each other at the infinity.

\para

More precisely, we say that  a  curve $\overline{\cal C}$ {\em
approaches} ${\cal C}$ at its infinity branch $B$ if the distance between
$\overline{\cal C}$ and $B$ approaches zero as they tend to infinity. We provide some results that
characterize whether two  plane algebraic curves are approaching.

\para

Using these results, we present a method to compare the
asymptotic behavior  of two curves (i.e. the behavior of two curves at the infinity). In particular, we prove that if two plane algebraic curves have
the same {\it asymptotic behavior}, the Hausdorff distance between
them is finite. As a consequence of the  results obtained in this paper, in \cite{paper2}, we present an algorithm for computing all the {\it generalized
asymptotes} of a real plane algebraic curve $\cal C$ defined
implicitly. The algorithm is based on the notion of perfect curve that, intuitively speaking, defines a curve of degree $d$ that cannot be approached by
any curve of degree less than $d$.

\para

The structure of the paper is as follows: In Section 2, we present the terminology that will be used throughout this paper as well as
some previous results. In Section 3,  the notion of {\it infinity
branch}  is introduced and
some important properties are proved.   In Section 4, we provide
the notions of {\it convergent branches} and {\it approaching
curves}. In addition, we develop some results that characterize whether
two plane algebraic curves approach each other. The results
presented in this section will be used in Section 5, where  an algorithm
to compare the asymptotic behavior of two algebraic plane curves is developed.  In
addition, we prove that if two plane   curves have the same
asymptotic behavior,   the Hausdorff distance between them is
finite.

\section{Preliminaries and Terminology}

In this section, we present some notions and terminology that will
be used throughout the paper. In particular, we need some previous
results concerning local parametrizations and Puiseux series.  For further details see \cite{Duval89}, Section 2.5 in \cite{SWP},
\cite{Stad00},  and  Chapter 4 (Section 2) in \cite{Walker}.

\para

We denote by ${\Bbb C}[[t]]$ the domain of {\em formal power series}
in the indeterminate $t$ with coefficients in the field ${\Bbb C}$,
i.e. the set of all sums of the form $\sum_{i=0}^\infty a_it^i$,
 $a_i \in {\Bbb C}$. The quotient field of ${\Bbb C}[[t]]$ is
called the field of {\em formal Laurent series}, and it is denoted by
${\Bbb C}((t))$. It is well known that every non-zero formal Laurent
series $A \in {\Bbb C}((t))$ can be written in the form
$A(t) = t^k\cdot(a_0+a_1t+a_2t^2+\cdots), {\rm\ where\ }
    a_0\not= 0 {\rm\ and\ } k\in \Bbb Z.$
  In addition, the field
${{\Bbb C}\ll t\gg}:= \bigcup_{n=1}^\infty {\Bbb C}((t^{1/n}))$ is
called the field of {\em  formal Puiseux series}. Note that Puiseux
series are power series with fractional exponents. In addition,
every Puiseux series, $\varphi$, has a bound
 for the denominators of exponents
with non-vanishing coefficients, which is known as
 {\em  the ramification index} of the series. We denote it as $\nu(\varphi)$ (see \cite{Duval89}).

\para

The {\em order} of a non-zero (Puiseux or Laurent) series $A$ is the
smallest exponent of a term with non-vanishing coefficient in $A$.
We denote  it  by $\ord(A)$. We let the order of 0 be $\infty$.

\para

In the following, we introduce the notion of projective local pa\-ra\-me\-tri\-zation
for a projective plane curve  (see Definition 2.69, and Lemma 2.70 in
\cite{SWP}).

\begin{definition}
Let ${\cal C}^* \subset {\Bbb P}^2({\Bbb C})$ be a projective plane curve defined by
the homogeneous polynomial $F(x,y,z) \in {\Bbb R}[x,y,z]$. Let
$A^*,B^*,C^*$ be series in ${\Bbb C}((t))$ such that: (i) $F(A^*(t):B^*(t):C^*(t))=0$ (where the three series converge), and (ii)  there is no $D \in {\Bbb C}((t))\setminus \{0\}$ such that
    $D\cdot(A^*,B^*,C^*) \in {\Bbb C}^3$. Then ${\cal P^*}=(A^*:B^*:C^*) \in {\Bbb P}^2({\Bbb C}((t))\,)$ is
called a {\em projective local parametrization} of $\cal C^*$. In
addition, one can always find such a parametrization having
$\min\{\ord(A^*),\ord(B^*),\ord(C^*)\}=0$, and the point ${\cal P^*}(0)\in
\cal C^*$ is called the {\em center} of $\cal P^*$.
\end{definition}

\para

\noindent
For an affine plane curve, the above notion can be stated as follows:

\begin{definition}
Let ${\cal C}$ be a real plane algebraic curve over ${\Bbb C}$   defined implicitly by
the irreducible polynomial $f(x,y) \in {\Bbb R}[x,y]$. Let $A,B$ be series in
${\Bbb C}((t))$ such that: (i) $f(A(t),B(t))=0$ (where both series converge), and
(ii) not both, $A$ and $B$, are constants.
Then ${\cal P}=(A,B)$ is called an {\em (affine) local
parametrization} of $\cal C$. Moreover, if $\ord(A),\ord(B)\geq 0$,
  the point ${\cal P}(0)=(a,b)\in {\cal C}$ is called the {\em
center} of $\cal P$.
\end{definition}

In the following, we deal with affine curves. The results and notions presented can be adapted for projective curves in an obvious way.

\para

Two local parametrizations, ${\cal P}_1$ and ${\cal P}_2$, of an
algebraic plane curve $\cal C$ are called {\em equivalent} if there exists
$R\in {\Bbb C}[[t]]$, with $\ord(R)=1$, such that ${\cal P}_1 ={\cal
P}_2(R)$. It can be proved that this equivalence of local
parametrizations is actually an equivalence relation.

\para

%Under these conditions, one has the following theorem (see  Chapter 4, Subsection 2.1, in \cite{Walker}):

%\begin{theorem}\label{T-localpara}
%  In a suitable affine coordinate system, any given local
%para\-metrization is equivalent to one of the form $\bigl(t^n,\
%\psi(t)\bigr)$, where
%$\psi(t)=a_1t^{n_1}+a_2t^{n_2}+a_3t^{n_3} + \cdots$,\,$n,n_i\in {\Bbb N},\,\,i=1,\ldots$, and $0<n_1<n_2<\cdots$.
%\end{theorem}

If a local parametrization ${\cal P}$, or one equivalent,
satisfies that ${\cal P}(t) = {\cal P}'(t^k)$ for some
parametrization ${\cal P}'$ and for some natural number $k>1$, then
${\cal P}$ is said to be {\em reducible}. Otherwise, ${\cal P}(t)$
is said to be {\em irreducible}. Under these conditions, we introduce
the notion of {\em place} as follows.

\begin{definition} An equivalence class of irreducible
local parametrizations of the algebraic plane curve $\cal C$ is called a
{\em place} of $\cal C$. The common center of the local
parametrizations (if it exists) is the {\em center} of the place.
\end{definition}

\para

In the following definition, we introduce the notion of {\em branch} of a plane curve.

\begin{definition} Given a local parametrization $(X,Y)$ of a plane curve $\cal C$,
the set of all points $(X(t),Y(t))$ obtained by allowing $t$ to vary
within some neighborhood of $0$ where $X(t)$ and $Y(t)$ converge is
called a {\em branch} of $\cal C$.
\end{definition}

It can be shown that two equivalent local parametrizations
provide the same branch. Therefore, one obtains a branch for each place
of a given algebraic plane curve.

\para

One may prove that the center of a local parametrization of $\cal C$
is a point on $\cal C$. Conversely,   from the following theorems,
we also obtain that every point on $\cal C$ is the center of at
least one place of $\cal C$  (see Theorems 2.77 and 2.78 in
\cite{SWP}).

\begin{theorem}\label{T-Puiseux} {\bf (Puiseux's Theorem)}
 The field $K\ll x\gg$ is algebraically closed.
\end{theorem}

A proof of Puiseux's Theorem can be given constructively by the
Newton Polygon Method (see e.g. Section 2.5 in \cite{SWP}). This method solves the construction of solutions of non-constant
univariate polynomial equations   over $K\ll x\gg$.

\begin{theorem} \label{T-puiseux-place}
Let  $\cal C$ be a plane curve defined by $f(x,y)\in {\Bbb R}[x,y]$.
To each root $Y(x) \in {\Bbb C}\ll x\gg$ of $f(x,y)=0$ with
$\ord(Y)>0$ there corresponds a unique place of $\cal C$ with center
at the origin. Conversely, to each place $(X(t),Y(t))$ of $\cal C$
with center at the origin there correspond $\ord(X)$ roots of
$f(x,y)=0$, each of order greater than zero.
\end{theorem}

\para

If $Y(x)$ is a Puiseux series solving $f(x,y)=0$, $\ord(Y)>0$, and
$n$ is the least integer for which $Y(x) \in {\Bbb C}((x^{1\over
n}))$  (i.e., $\nu(Y)=n$),
then we set $x^{1\over n}=t$, and $(t^n, Y(t^n))$ is a local
parametrization with center at the origin. The solutions of $f(x,y)$
of order 0 correspond to places with center on the $y$-axis but
different from the origin, and the solutions of negative order
correspond to places at infinity (places with center at an infinity
point).

\para

Note that several different Puiseux series may correspond to
equivalent local parametrizations, and then these series provide a
unique place.  More precisely, let $Y(x)=\sum_{i\geq r}a_ix^{i/n}$ be a Puiseux
series with ramification index $\nu(Y)=n$. The series $\sigma_{\epsilon}(Y)$, $\epsilon^n=1$, are called
the {\em conjugates} of $Y$, where
$$\sigma_{\epsilon}(Y)=\sum_{i\geq r}\epsilon^ia_ix^{i/n}.$$
The set of all (distinct) conjugates of $Y$ is called the  {\em
conjugacy class} of $Y$. The number of different conjugates of $Y$
is $\nu(Y)$. Two Puiseux series provide the same place if they
belong to the same conjugacy class (see \cite{Duval89} and
\cite{verger}).

%So we have the problem of identifying equivalent Puiseux series.
%This question is dealt by some authors as for instance in
%\cite{Duval89}.
% it is  developed a refinement of the classical Newton polygon method  which allows to detect this kind of parameter substitutions

\section{Infinity Branches}\label{sec-Infinity branches}

In this section, we introduce the notion of {\it infinity branch}
(see Definition \ref{D-infinitybranch}), and we obtain some  properties concerning to these algebraic entities.

\para

For this purpose, we consider an algebraic affine plane curve $\cal C$  over
$\mathbb{C}$, defined implicitly by the irreducible polynomial $f(x,y) \in {\Bbb R}[x,y]$.
Let ${\cal C}^*$  be its corresponding projective curve  defined by
the homogeneous polynomial $F(x,y,z) \in {\Bbb R}[x,y,z]$. Furthermore, let $P=(1:m:0),\,m\in {\Bbb C},$ be an infinity point of  ${\cal C}^*$, and we consider the  curve   defined implicitly by the
polynomial $g(y,z)=F(1:y:z)$. Observe that $g(p)=0,$ where $p=(m,0)$.

\para

By applying Theorem  \ref{T-Puiseux}, we compute the
series expansion for the solutions of $g(y,z)=0$. There exist
exactly $\deg_{Y}(g)$ solutions given by different Puiseux series
that can be grouped into conjugacy classes. Let one of these
solutions be given by the following Puiseux series:
$$\varphi(z)=m+a_1z^{N_1/N}+a_2z^{N_2/N}+a_3z^{N_3/N} +
\cdots\in{\Bbb C}\ll z\gg,\quad a_i\not=0,\, \forall i\in {\Bbb
N},$$
 where $\nu(\varphi)=N\in {\Bbb
N}$, $N_i\in {\Bbb N},\,\,i=1,\ldots$,  and $0<N_1<N_2<\cdots$. We
have that $g(\varphi(z), z)=0$ in some neighborhood of $z=0$ where
$\varphi(z)$ converges. Then, there exists some $M \in {\Bbb R}^+$
such that
$$F(1:\varphi(t):t)=g(\varphi(t), t)=0,\quad \mbox{for\, $t\in {\Bbb C}$\, and
$|t|<M$},$$
which implies that
$F(t^{-1}:t^{-1}\varphi(t):1)=f(t^{-1},t^{-1}\varphi(t))=0$, for
$t\in {\Bbb C}$  and  $0<|t|<M$. We set $t^{-1}=z$, and   we obtain
that
$$f(z,r(z))=0,\quad \mbox{$z\in {\Bbb C}$\, and
$|z|>M^{-1}$,\qquad where}$$
$$r(z)=z\varphi(z^{-1})=mz+a_1z^{1-N_1/N}+a_2z^{1-N_2/N}+a_3z^{1-N_3/N} + \cdots,\quad a_i\not=0,\, \forall i\in {\Bbb N} $$
$N,N_i\in {\Bbb N},\,\,i=1,\ldots$, and $0<N_1<N_2<\cdots$.

\para

  Since $\nu(\varphi)=N$, we get that there are $N$ different series in its conjugacy
class. Let  $\varphi_1,\ldots,\varphi_N$ be these series,  and
\begin{equation}\label{Eq-conjugates}r_i(z)=z\varphi_i(z^{-1})=mz+a_1c_i^{N_1}z^{1-N_1/N}+a_2c_i^{N_2}z^{1-N_2/N}+a_3c_i^{N_3}z^{1-N_3/N}
+ \cdots\end{equation} where $c_1,\ldots,c_N$ are the $N$ complex
roots of $x^N=1$. Now we are ready to introduce the notion of
infinity branch.

\para

\begin{definition}\label{D-infinitybranch} The set
$\displaystyle B=\bigcup_{i=1}^NL_i$ where
$$L_i=\{(z,r_i(z))\in {\Bbb
C}^2: \,z\in {\Bbb C},\,|z|>M_i\}$$ is called an {\em infinity
branch} of the affine plane curve $\cal C$. The subsets
$L_1,\ldots,L_N$ are called the {\em leaves} of the infinity branch
$B$.
\end{definition}

\para

  \begin{remark} \begin{enumerate}
\item We observe that an infinity branch  is uniquely determined  from one leaf, up to conjugation. That is, if
$\displaystyle B=\bigcup_{i=1}^NL_i$, where $L_i=\{(z,r_i(z))\in
{\Bbb C}^2: \,z\in {\Bbb C},\,|z|>M_i\}$, and
$$ r_i(z)=z\varphi_i(z^{-1})=mz+a_1z^{1-N_1/N}+a_2z^{1-N_2/N}+a_3z^{1-N_3/N} + \cdots$$
then $r_j=r_i,\,j=1,\ldots,N$, up to conjugation; i.e.
\[r_j(z)=z\varphi_j(z^{-1})=mz+a_1c_j^{N_1}z^{1-N_1/N}+a_2c_j^{N_2}z^{1-N_2/N}+a_3c_j^{N_3}z^{1-N_3/N}
+ \cdots\] where $c_j^N=1,\,\,j=1,\ldots,N$, and
$N,N_i\in\mathbb{N}$.
\item Let  $M:=\max\{M_1,\ldots,M_N\}$.  In the following, we consider $L_i=\{(z,r_i(z))\in {\Bbb
C}^2: \,z\in {\Bbb C},\,|z|>M\}$.\end{enumerate}\end{remark}

\para

Let $\varphi(z)=m+a_1z^{N_1/N}+a_2z^{N_2/N}+a_3z^{N_3/N} + \cdots$
be a series expansion for a solution  of $g(y,z)=0$. We consider
$\psi(t):=\varphi(t^N),$ and we observe that $(1:\psi(t):t^N)$ is a
local projective parametrization, with center at $P$, of the
projective curve ${\cal C}^*$.

\para

Thus, from $\psi_i(t):=\varphi_i(t^N),\,\,i=1,\ldots,N$ ($\varphi_i$
are the $N$ different series in the conjugacy class of $\varphi$),
we obtain $N$ equivalent local projective parametrizations,
$(1:\psi_i(t):t^N)$ (note that they are equivalent since
$\varphi_1,\ldots\varphi_N$ belong to the same conjugacy class).
Therefore, the leaves of $B$ are all associated to a unique infinity
place.

\para

Conversely, from a given  infinity place defined  by a local
projective pa\-ra\-me\-tri\-zation $(1:\psi(t):t^N)$ (see Theorem
2.5.3 in \cite{SWP}), we obtain $N$ Puiseux series,
$\varphi_j(t)=\psi(c_j{t}^{1/N})$, $c_j^N=1$, that provide different
expressions $r_j(z)=z\varphi_j(z^{-1}),\,j=1,\ldots,N$. Hence, the
infinity branch $B$ is defined by the leaves $L_j=\{(z,r_j(z))\in
{\Bbb C}^2: \,z\in {\Bbb C},\,|z|>M\},\,\,j=1,\ldots,N.$

\para

From the above discussion, we deduce that there exists a
one-to-one relation between infinity places and infinity branches.
In addition, we can say that each infinity branch is associated to
a unique infinity point given by the center of the corresponding infinity place. Reciprocally,
taking into account the above construction, we get that every infinity point has associated, at least, one
infinity branch. Hence, every algebraic plane curve has, at least, one infinity
branch. Furthermore,  every algebraic plane curve has a finite number of branches.

\para

Observe that the above construction can be applied to any infinity point of the form $(a:b:0),\,a\not=0$. In the following, we assume that $a=0$; that is, we take  the infinity point $P=(0:1:0)$. In this case,
 we consider the  curve  defined implicitly by the
polynomial $h(x,z)=F(x:1:z)$. Observe that $h(p)=0,$ where $p=(0,0)$.
In this situation, we get that there exists $M \in {\Bbb R}^+$ such that
$$F(\varphi(t):1:t)=h(\varphi(t), t)=0,\quad \mbox{for\, $t\in {\Bbb C}$\, and
$|t|<M$,\,\,\quad where}$$
$$\varphi(z)=a_1z^{N_1/N}+a_2z^{N_2/N}+a_3z^{N_3/N} + \cdots\in{\Bbb C}\ll z\gg,\quad a_i\not=0,\, \forall i\in {\Bbb N}$$\,
$N,N_i\in {\Bbb N},\,\,i=1,\ldots$, and $0<N_1<N_2<\cdots$,  is a
series expansion for a solution  of $h(x,z)=0$. We set $z=t^{-1}$,
and we get that
$$f(r(z),z)=0,\quad \mbox{$z\in {\Bbb C}$\, and
$|z|>M^{-1}$,\qquad where}$$
$$r(z)=z\varphi(z^{-1})=a_1z^{1-N_1/N}+a_2z^{1-N_2/N}+a_3z^{1-N_3/N} + \cdots,\quad a_i\not=0,\, \forall i\in {\Bbb N}$$
$N,N_i\in {\Bbb N},\,\,i=1,\ldots$, and $0<N_1<N_2<\cdots$.

\para

\noindent Thus, we obtain an infinity branch $\displaystyle
B=\bigcup_{i=1}^NL_i$ whose leaves have the form:
$$L_i=\{(r_i(z),z)\in {\Bbb
C}^2: \,z\in {\Bbb C},\,|z|>M\}.$$
Observe that we may apply this   construction to any infinity point of the form $(a:b:0),\,b\not=0$.

\para

\noindent
These two approaches lead  us to consider two types of infinity
branches.

\begin{definition}  Let  ${\cal C}$ be an affine plane curve over ${\Bbb C}$ defined by an irreducible
polynomial $f(x,y) \in {\Bbb R}[x,y]$.
\begin{itemize}

\item An {\em
infinity branch of ${\cal C}$ of type 1}  associated to the infinity
point $P=(1:m:0),\,m\in {\Bbb C}$, is  a set $\displaystyle
B=\bigcup_{i=1}^NL_i$, where $L_i=\{(z,r_i(z))\in {\Bbb C}^2: \,z\in
{\Bbb C},\,|z|>M\}$,\, $i=1,\ldots,N$, $M\in {\Bbb R}^+$, and
$r_1,\ldots,r_N$ are the conjugates of
$$r(z)=mz+a_1z^{1-N_1/N}+a_2z^{1-N_2/N}+a_3z^{1-N_3/N} + \cdots,\quad a_i\not=0,\, \forall i\in {\Bbb N}$$
$N,N_i\in {\Bbb N},\,\,i=1,\ldots$, and $0<N_1<N_2<\cdots$.

\item An  {\em infinity branch of
${\cal C}$ of type 2}  associated to the infinity point
$P=(m:1:0),\,m\in {\Bbb C}$,  is  a set $\displaystyle
B=\bigcup_{i=1}^NL_i$, where $L_i=\{(r_i(z),z)\in {\Bbb C}^2: \,z\in
{\Bbb C},\,|z|>M\}$,\, $i=1,\ldots,N$, $M\in {\Bbb R}^+$, and
$r_1,\ldots,r_N$ are the conjugates of
$$r(z)=mz+a_1z^{1-N_1/N}+a_2z^{1-N_2/N}+a_3z^{1-N_3/N} + \cdots,\quad a_i\not=0,\, \forall i\in {\Bbb N}$$
$N,N_i\in {\Bbb N},\,\,i=1,\ldots$, and $0<N_1<N_2<\cdots$.

\end{itemize}
\end{definition}

\para

\begin{remark} \begin{enumerate}
\item In the following, we assume w.l.o.g that  the given algebraic plane curve $\cal C$ only has  type 1 infinity branches; that is,  all the
infinity points are of the form $(1:m:0)$, $m\in {\Bbb C}$. Otherwise, we may consider a linear change of coordinates.
\item  By abuse of notation, we will say that $N$ is  {\em  the ramification index} of the branch $B$, and we will write it as
$\nu(B)=N$. Note that $B$ has $\nu(B)$  leaves.
\end{enumerate}
\end{remark}

In the following example, we compute the infinity branches for a given plane curve.

\begin{example}\label{Ej-ramas infinitas}
Let  ${\cal C}$ be the plane curve defined implicitly by the irreducible
polynomial
$$f(x,y)=y^5-4y^4x+4y^3x^2+2y^2x-y^2x^2+2yx^2+2yx^3+x+x^2\in {\Bbb R}[x,y].$$
The corresponding   projective
curve ${\cal C}^*$ is defined by  $F(x:y:z)=$
$$y^5-4y^4x+4y^3x^2+2y^2z^2x-zy^2x^2+2z^2yx^2+2yzx^3+z^4x+z^3x^2\in {\Bbb R}[x,y,z].$$
Note that $P=(1:0:0)$ is an infinity point of   ${\cal C}^*$. Let us compute the infinity branches associated to $P$. For this purpose,   we consider the  curve defined implicitly by the
polynomial $g(y,z)=F(1:y:z)$, and we observe that $g(p)=0,$ where $p=(0,0)$.

\para

\noindent
We compute the   series
expansion for the solutions of $g(y,z)=0$. For this purpose, we use for instance  the {\sf algcurves} package included in  the computer algebra system {\sf Maple}. We get that:
$$\varphi_1(z)=-1/2z^2+1/8z^4-1/8z^5+1/16z^6+1/16z^7+\cdots\in{\Bbb C}\ll z\gg,\quad \mbox{and}$$
$$\varphi_2(z)=\frac{-(-2z)^{1/2}}{2}-\frac{z}{8}  +\frac{27}{256} (-2 z)^{3/2}-\frac{7}{32} z^2+\frac{4057}{65536} (-2 z)^{5/2}+\cdots\in{\Bbb C}\ll z\gg.$$

\para

That is,  $g(\varphi_j(z), z)=0,\,j =1,2$ (see e.g. Section 2.5 in \cite{SWP}).
Note that   $\nu(\varphi_1)=1$, which implies that we only have one
Puiseux series in the conjugacy
class  of $\varphi_1$.  However,  $\nu(\varphi_2)=2$ and then, we have the following conjugate
Puiseux series in the conjugacy
class  of $\varphi_2$:
$$\varphi_{2,1}(z)=\frac{-(-2z)^{1/2}}{2}-\frac{z}{8}  +\frac{27}{256} (-2 z)^{3/2}-\frac{7}{32} z^2+\frac{4057}{65536} (-2
z)^{5/2}+\cdots$$
$$\varphi_{2,2}(z)=\frac{+(-2z)^{1/2}}{2}-\frac{z}{8}  -\frac{27}{256} (-2 z)^{3/2}-\frac{7}{32} z^2-\frac{4057}{65536} (-2
z)^{5/2}+\cdots$$
%Thus, there exists $M \in {\Bbb R}^+$ such that
%$$F(1:\varphi(t):t)=g(\varphi(t), t)=0,\quad \mbox{for\, $t\in {\Bbb C}$\, and
%$|t|<M$}$$  which implies that $F(t^{-1}:t^{-1}\varphi(t):1)=f(t^{-1},t^{-1}\varphi(t))=0$, for
%$t\in {\Bbb C}$  and  $0<|t|<M$.
%\para
%\noindent
%Now, we set $t^{-1}=z$, and   we obtain that
%$$f(z,r(z))=F(z:r(z):1)=0,\quad \mbox{$z\in {\Bbb C}$\, and
%$|z|>M^{-1}$,}$$ where
%$$r(z)=z\varphi(z^{-1})=mz+a_1z^{1-n_1/n}+a_2z^{1-n_2/n}+a_3z^{1-n_3/n} + \cdots, $$
%$n,n_i\in {\Bbb N},\,\,i=1,\ldots$, and $0<n_1<n_2<\cdots$.
Thus, we obtain two infinity branches:
$$B_1=L_1=\{(z,r_1(z))\in {\Bbb C}^2: \,z\in {\Bbb C},\,|z|>M\},\qquad \mbox{where}$$
$$r_1(z)=z\varphi_1(z^{-1})=-1/(2z)+1/(8z^3)-1/(8z^4)+1/(16z^5)+1/(16z^6)+\cdots$$
and $B_2=L_{2,1}\cup L_{2,2}$, where $L_{2,i}=\{(z,r_{2,i}(z))\in {\Bbb C}^2:
\,z\in {\Bbb C},\,|z|>M\}$, $i=1,2$ and
$$r_{2,1}(z)=z\varphi_{2,1}(z^{-1})=-\frac{(-2z)^{1/2}}{2}-\frac{1}{8}+\frac{ 27(-2z)^{-1/2}}{64}-\frac{7z^{-1}}{32}+\frac{ 4057(-2z)^{-3/2}}{4096}+\cdots.$$
$$r_{2,2}(z)=z\varphi_{2,2}(z^{-1})=+\frac{(-2z)^{1/2}}{2}-\frac{1}{8}-\frac{ 27(-2z)^{-1/2}}{64}-\frac{7z^{-1}}{32}-\frac{
4057(-2z)^{-3/2}}{4096}+\cdots.$$

\para

 In Figure \ref{F-ramas infinitas2}, we plot the curve
${\cal C}$ and some points of the infinity branches $B_1$ and $B_2$ associated to $P$.

\begin{figure}[h]
$$
\begin{array}{cc}
\psfig{figure=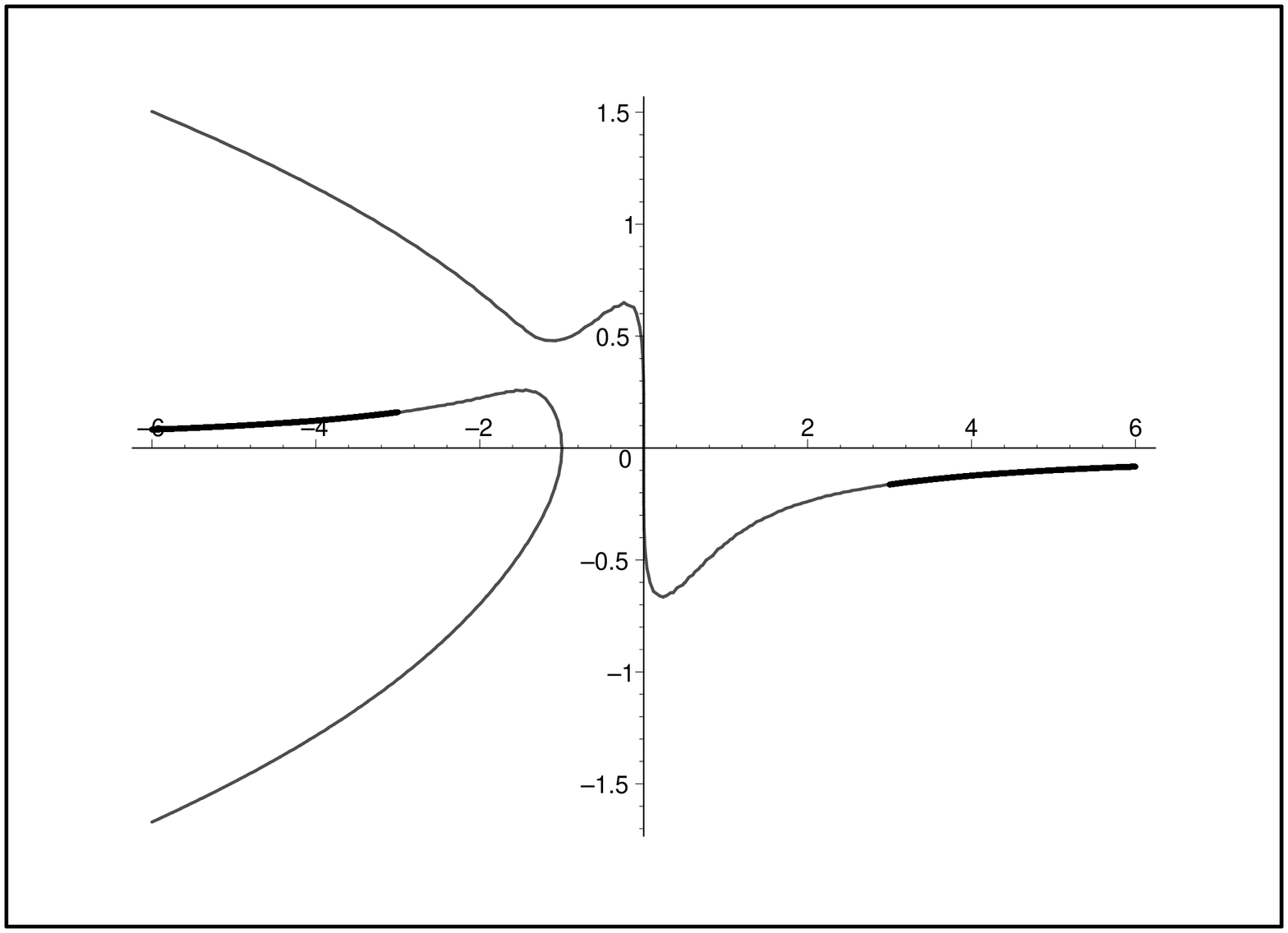,width=4.5cm,height=4.5cm,angle=270} &
\psfig{figure=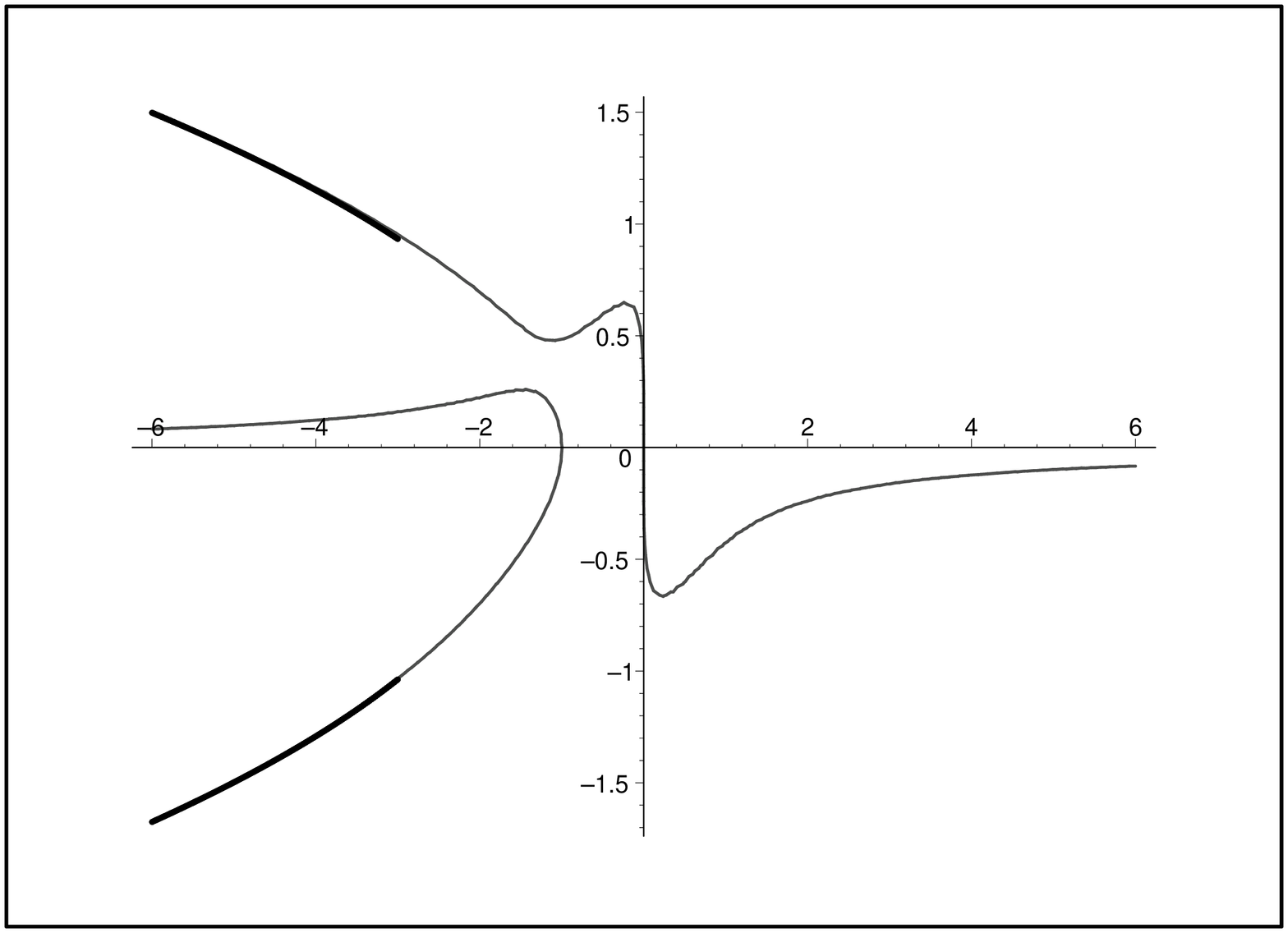,width=4.5cm,height=4.5cm,angle=270}
\end{array}
$$ \caption{Infinity branches $B_1$ (left), and $B_2$
(right).}\label{F-ramas infinitas2}
\end{figure}
\end{example}

 \para

In the following, we prove that  any point of the curve with
sufficiently large coordinates  belongs to some infinity branch. For this purpose, we recall the reader   that if $h$ is a complex-valued function of a complex
variable,  $h: {\Bbb C} \rightarrow {\Bbb C}$,  we say that the limit of $h(z)$ as $z$ approaches $\infty$ is $L$, written
$\displaystyle\lim_{z\rightarrow\infty}h(z)=L$, if whenever  $\{z_n\}_{n\in {\Bbb N}}$  is a sequence of points with $\displaystyle\lim_{n\rightarrow\infty} z_n =\infty$,  it holds that
 $\displaystyle\lim_{n\rightarrow\infty}h(z_n)=L$ (see e.g. \cite{Ahlfors} or \cite{conway}).

 \para

\begin{lemma}\label{L-BolaMaxima} Let $\cal C$ be an algebraic plane curve.
There exists $K\in\mathbb{R}^+$ such that for every
$p=(a,b)\in\mathcal{C}$ with $|a|>K$, it holds that $p\in B_p$, where $B_p$ is an infinity branch of ${\cal C}$.
\end{lemma}

\noindent\textbf{Proof:} Let us assume that the lemma does not hold, and we
consider a sequence $\{K_n\}_{n\in {\Bbb N}}\subset \mathbb{R}^+$ such that
$\lim_{n\rightarrow\infty}K_n=\infty$. Then, for every $n\in {\Bbb N}$ there   exists a point
$p_n=(a_n,b_n)\in\mathcal{C}$ such that $|a_n|>K_n$, and $p_n$ does not
belong to any infinity branch of $\cal C$.

\para

Let $P_n=(a_n:b_n:1)$. Since  $F(P_n)=f(p_n)=0$, then
$\lim_{n\rightarrow\infty}F(P_n)=0$. Thus, we
distinguish two different cases:

\begin{enumerate}
\item[a)] If there   exists a monotone subsequence $\{b_{n_l}/a_{n_l}\}_{l\in {\Bbb N}}$ that is not bounded, we have that
 $\lim_{l\rightarrow\infty}b_{n_l}/a_{n_l}=\infty$, and then
  $\lim_{l\rightarrow\infty}a_{n_l}/b_{n_l}=0$. Hence,
$$\lim_{l\rightarrow\infty}F(Q_{n_l})=F(0:1:0)=0,\qquad Q_{n_l}=(a_{n_l}/b_{n_l}:1:1/b_{n_l})$$
which implies that
$P=(0:1:0)$ is an infinity point of ${\cal C}^*$.
\item[b)]  If there   exists  a monotone subsequence $\{b_{n_l}/a_{n_l}\}_{l\in {\Bbb N}}$ that is bounded, we have that  $\lim_{l\rightarrow\infty}b_{n_l}/a_{n_l}=m$. Thus,
$$\lim_{l\rightarrow\infty}F(Q_{n_l})=F(1:m:0)=0,\qquad Q_{n_l}=(1:b_{n_l}/a_{n_l}:1/a_{n_l})$$
which implies that  $P=(1:m:0)$ is an infinity point of ${\cal C}^*$.
\end{enumerate}

From both situations, we conclude that there exist a sequence
$\{Q_n\}_{n\in {\Bbb N}}$  that approaches to an infinity point $P$
as $n$  tends to infinity; that is, there exists $M\in {\Bbb R}^+$
such that $\|Q_n-P\|\leq \epsilon$, for $n\geq M$. Thus, %from the Implicit Function Theorem,
 we deduce that  $\{Q_n\}_{n\in {\Bbb
N},\,n\geq M}$   can be obtained by a place centered at $P$. Hence,
$p_n$ belongs to some infinity branch of ${\cal C}$, which
contradicts the hypothesis.
 \hfill $\Box$

 \para

\begin{remark}\label{R-BolaMaxima} Reasoning similarly as in Lemma \ref{L-BolaMaxima}, one has that there exists $K\in\mathbb{R}^+$ such that for every
$p=(a,b)\in\mathcal{C}$ with $|b|>K$, it holds that $p\in B_p$, where   $B_p$ is an infinity branch of ${\cal C}$.
\end{remark}

\section{Convergent Branches and \\Approaching Curves}

In this section, we introduce the notions of convergent branches and
approaching curves. Intuitively speaking, two infinity branches
converge if they get closer  as they tend to infinity. This
concept will allow us to analyze whether two curves approach each
other at the infinity.
\para

The results presented in this section will be used in Section 5, where  a method to compare the
asymptotic behavior of two curves is developed.

\para

 \begin{definition} Given two leaves, $L=\{(z,r(z))\in
{\Bbb C}^2:\,z\in {\Bbb C},\,|z|>M\}$ and
$\overline{L}=\{(z,\overline{r}(z))\in {\Bbb C}^2:\,z\in {\Bbb
C},\,|z|>\overline{M}\}$, we say that they are convergent if
$\lim_{z\rightarrow\infty} (\overline{r}(z)-r(z))=0.$
\end{definition}

%In the following lemma, we characterize the convergence of two
%\textcolor{red}{leaves}.

 \para

\begin{lemma}\label{L-DistVertical} Two  leaves $L=\{(z,r(z))\in {\Bbb C}^2:\,z\in {\Bbb C},\,|z|>M\}$ and
$\overline{L}=\{(z,\overline{r}(z))\in {\Bbb C}^2:\,z\in {\Bbb
C},\,|z|>\overline{M}\}$ are convergent if and only if the terms
with non negative exponent in the series $r(z)$ and
$\overline{r}(z)$ are the same.
\end{lemma}

\noindent\textbf{Proof:}  Let
$$r(z)=mz+a_1z^{\frac{N-N_1}{N}}+a_2z^{\frac{N-N_2}{N}}+\cdots ,\,\, N, N_i\in {\Bbb N},\, 0<N_1<N_2<\cdots,\, \, a_i\not=0$$ and
$$ \overline{r}(z)=\overline{m}z+b_1z^{\frac{\overline{N}-\overline{N}_1}{\overline{N}}}+b_2z^{\frac{\overline{N}-\overline{N}_2}{\overline{N}}}
+\cdots,\,\, \overline{N}, \overline{N}_i\in {\Bbb N},\,
0<\overline{N}_1<\overline{N}_2<\cdots,\,\,b_i\not=0. $$ Then,
$$r(z)-\overline{r}(z)=mz-\overline{m}z+a_1z^{\frac{N-N_1}{N}}-b_1z^{\frac{\overline{N}-\overline{N}_1}{\overline{N}}}
+a_2z^{\frac{N-N_2}{N}}-b_2z^{\frac{\overline{N}-\overline{N}_2}{\overline{N}}}+\cdots.$$
Note that $\lim_{z\rightarrow\infty}(r(z)-\overline{r}(z))=0$ if and
only if $r(z)-\overline{r}(z)$ has no terms with non negative
exponent. This situation holds if the terms with non negative
exponent in both series, $r(z)$ and $\overline{r}(z)$, are the same.
\hfill $\Box$

\para

\begin{remark}\label{R-DistVertical} \begin{enumerate}
\item  From Lemma \ref{L-DistVertical}, we deduce that $m=\overline{m}$ and then,  $L$ and $\overline{L}$
are associated to the same infinity point.
\item Note that the number of terms with positive exponent in both series is finite.
\end{enumerate}
\end{remark}

\para

 \begin{definition}\label{D-distance0} Two
infinity branches, $B$ and $\overline{B}$, are convergent if there
exist two convergent leaves $L\subset B$ and $\overline{L}\subset
\overline{B}$.
\end{definition}

  \begin{remark}\label{R-common-inf-point} From Remark
\ref{R-DistVertical}, statement 1, we get that two convergent infinity branches
are associated to the same infinity point.
% Indeed: since
%$\lim_{z\rightarrow\infty}(\overline{r}(z)-r(z))=0$, we have that
%$\lim_{z\rightarrow\infty}\left(\frac{\overline{r}(z)}{z}-\frac{r(z)}{z}\right)=0$
%which implies that
%$\lim_{z\rightarrow\infty}\frac{\overline{r}(z)}{z}=\lim_{z\rightarrow\infty}\frac{r(z)}{z}.$
\end{remark}

\begin{proposition}\label{P-convergent-branches}
Two infinity branches $B$ and $\overline{B}$ are convergent if and
only if for each leaf $L\subset B$ there exists a leaf $\overline{L}\subset
\overline{B}$ convergent with $L$, and reciprocally.
\end{proposition}

\noindent\textbf{Proof:} Let $B$ and $\overline{B}$ be two convergent
infinity branches, and let us prove that for any $L_i\subset B$
there exists $\overline{L}_j\subset \overline{B}$ convergent with
$L_i$ (using Definition \ref{D-distance0}, we clearly have the reciprocal). From Definition \ref{D-distance0},
there exist two leaves $L=\{(z,r(z))\in {\Bbb C}^2:\,z\in {\Bbb
C},\,|z|>M\}\subset B$, and $\overline{L}=\{(z,\overline{r}(z))\in
{\Bbb C}^2:\,z\in {\Bbb C},\,|z|>\overline{M}\}\subset \overline{B}$
convergent. Let
$$r(z)=z\varphi(z^{-1})=mz+u_1z^{1-\frac{N_1}{N}}+\cdots+u_kz^{1-\frac{N_k}{N}}+u_{k+1}z^{1-\frac{N_{k+1}}{N}}+\cdots,\, \, u_i\not=0,$$
$$\overline{r}(z)=z\overline{\varphi}(z^{-1})=mz+\overline{u}_1z^{1-\frac{\overline{N}_1}{\overline{N}}}
+\cdots+\overline{u}_kz^{1-\frac{\overline{N}_k}{\overline{N}}}
+\overline{u}_{k+1}z^{1-\frac{\overline{N}_{k+1}}{\overline{N}}}+\cdots,\,
\, \overline{u}_i\not=0,$$ where $\nu(B)=N$,
$\nu(\overline{B})=\overline{N}$, $N_k\leq N<N_{k+1}$ and
$\overline{N}_k\leq \overline{N}<\overline{N}_{k+1}$.

\para

From Lemma \ref{L-DistVertical}, we deduce that the terms with non
negative exponent in $r$ and $\overline{r}$ must coincide. Thus,
$u_l=\overline{u}_l=a_l$, for $l=1,\ldots,k$, and
$$r(z)=mz+a_1z^{1-\frac{n_1}{n}}+\cdots+a_kz^{1-\frac{n_k}{n}}+u_{k+1}z^{1-\frac{N_{k+1}}{N}}+\cdots,\, \, a_i,u_i\not=0$$
$$\overline{r}(z)=mz+a_1z^{1-\frac{n_1}{n}}+\cdots+a_kz^{1-\frac{n_k}{n}}+\overline{u}_{k+1}z^{1-\frac{\overline{N}_{k+1}}{\overline{N}}}
+\cdots,\, \, a_i,\overline{u}_i\not=0,$$ where $n, n_i\in {\Bbb
N},\, 0<n_1<\cdots<n_k<n$,
$N_{k+1}>N$,$\overline{N}_{k+1}>\overline{N}$. Observe that  we have
simplified the non negative exponents such that
$\gcd(n,n_1,\ldots,n_k)=1$. That is,
 for $l=1,\ldots,k$, there are $b,\overline{b}\in\mathbb{N}$ such that   $N_l=bn_l$,
$N=bn$, $\overline{N}_l=\overline{b}n_l$, and
$\overline{N}=\overline{b}n$.

\para

Under these conditions, we observe that the different leaves of $B$
and $\overline{B}$ are obtained by conjugation on $r(z)$ and
$\overline{r}(z)$. That is (see equation (\ref{Eq-conjugates})),
$$r_i(z)=mz+u_1c_i^{N_1}z^{1-\frac{N_1}{N}}+\cdots+u_kc_i^{N_k}z^{1-\frac{N_k}{N}}+u_{k+1}c_i^{N_{k+1}}z^{1-\frac{N_{k+1}}{N}}+\cdots$$
$$\overline{r}_j(z)=mz+\overline{u}_1d_j^{\overline{N}_1}z^{1-\frac{\overline{N}_1}{\overline{N}}}+\cdots+
\overline{u}_kd_j^{\overline{N}_k}z^{1-\frac{\overline{N}_k}{\overline{N}}}
+\overline{u}_{k+1}d_j^{\overline{N}_{k+1}}z^{1-\frac{\overline{N}_{k+1}}{\overline{N}}}+\cdots,$$
where $c_1,\ldots,c_{N}$ are the $N$ complex roots of $x^{N}=1$, and
$d_1,\ldots,d_{\overline{N}}$ are the $\overline{N}$ complex roots
of $x^{\overline{N}}=1$.

\para

We simplify the exponents and, using that $u_l=\overline{u}_l=a_l,\,
l=1,\ldots,k$, we get that:
$$r_i(z)=mz+a_1c_i^{N_1}z^{1-\frac{n_1}{n}}+\cdots+a_kc_i^{N_k}z^{1-\frac{n_k}{n}}+u_{k+1}c_i^{N_{k+1}}z^{1-\frac{N_{k+1}}{N}}+\cdots$$
$$\overline{r}_j(z)=mz+a_1d_j^{\overline{N}_1}z^{1-\frac{n_1}{n}}+\cdots+a_kd_j^{\overline{N}_k}z^{1-\frac{n_k}{n}}
+\overline{u}_{k+1}d_j^{\overline{N}_{k+1}}z^{1-\frac{\overline{N}_{k+1}}{\overline{N}}}+\cdots.$$

\para

Hence, we only have to show that for each $i\in\{1,\ldots,N\}$ there
exists $j\in\{1,\ldots,\overline{N}\}$ such that
$c_i^{N_l}=d_j^{\overline{N}_l}$ for every $l=1,\ldots,k$. Indeed:
since $c_i,\,i=1,\ldots,N$ are the $N$ complex roots of $x^{N}=1$,
we have that  $c_i=e^{\frac{2(i-1)\pi I}{N}}$, where $I$ is the
imaginary unit. Taking into account that $N=bn$, we deduce that
$c_i^b=e^{\frac{2(i-1)\pi I}{n}}$,\,$i=1,\ldots,N$, and
$c_i^b=c_{i+(m-1)n}^b$ for each $i=1,\ldots,n$, and $m=1,\ldots,b$.
That is, $(c_i^{b})^n=1$,\,$i=1,\ldots,n$. Reasoning similarly, we
have that $d_j^{\overline{b}}=e^{\frac{2(j-1)\pi
I}{n}}$,\,$j=1,\ldots,\overline{N}$, and
$d_j^{\overline{b}}=d_{j+(m-1)n}^{\overline{b}}$ for each
$j=1,\ldots,n$, and $m=1,\ldots,\overline{b}$. That is,
$(d_j^{\overline{b}})^n=1$,\,$j=1,\ldots,n$.  Therefore,
$c_i^b=d_{i+(m-1)n}^{\overline{b}}$,\,$m=1,\ldots,\overline{b}$, and
using that $N_l=bn_l$ and
$\overline{N}_l=\overline{b}n_l$,\,$l=1,\ldots,k$, it follows that
$c_i^{N_l}=d_{j}^{\overline{N}_l}$,\,$j=i+(m-1)n,\,m=1,\ldots,\overline{b}$.
 \hfill $\Box$

%{\color{blue}{Hence, we only have to show that for each $i\in\{1,\ldots,N_1\}$
%there exists $j\in\{1,\ldots,N_2\}$ such that
%$c_i^{s_l}=d_j^{\overline{s}_l}$ for every $l=1,\ldots,n$. Indeed: since $c_1,\ldots,c_{N_1}$ are the $N_1$ complex
%roots of $x^{N_1}=1$, and  $N_1=bn$, we have that $(c_i^{b})^n=1$ which implies that  the
%set ${\cal A}:=\{c_i^{b}:i=1,\ldots,N_1\}$ contains the $n$ complex roots of
%$x^n=1$. In fact, ${\cal A}$ contains $n$ different values and each value appears $b$ times (note that for each $r$ such that $r^n=1$
%there are $b$ complex numbers such that $x^b=r$). Reasoning similarly, we deduce that  ${\cal B}:=\{d_j^{\overline{b}}:j=1,\ldots,N_2\}$ contains the
%$n$ different complex roots of $x^n=1$ (note that
%$d_j^{N_2}=(d_j^{\overline{b}})^n=1$), and each root appears $\overline{b}$ times.
%\para

%Thus, for each $c_i^{b}\in {\cal A}$ there are
%$\overline{b}$ values in ${\cal B}$ such that
%$c_i^{b}=d_j^{\overline{b}}$. Taking into account that
%$s_l=bn_l$ and $\overline{s}_l=\overline{b}n_l$, it follows that
%$c_i^{s_l}=d_j^{\overline{s}_l}$ for every $l=1,\ldots,n$. \hfill $\Box$}}

\para

\begin{remark}
Two convergent infinity branches may have different
 ramification indexes  i.e., they may have different number of
leaves. However, the value $n\in\mathbb{N}$ obtained by simplifying
the non negative exponents, is the same in both branches. We refer to it
the {\em degree of the infinity branch}. Observe that  the proof of
Proposition \ref{P-convergent-branches} implies that  two convergent
infinity branches have the same degree.\\  In order to
illustrate this remark, we consider the curves ${\cal C}$ and
$\overline{{\cal C}}$, defined  by the polynomials
$f(x,y)=y^4-2xy^2+x^2-y$, and $\overline{f}(x,y)=y^2-x$, respectively. ${\cal C}$ has only the infinity branch
$B=\bigcup_{i=1}^4 L_i$, where $L_i=\{(z,r_i(z))\in {\Bbb C}^2:
\,z\in {\Bbb C},\,|z|>M\}$,
$$r_i(z)=c_i^{2}
z^{1/2}+\frac{1}{2}c_i^{5}z^{-1/4}-\frac{1}{64}c_i^{11}z^{-7/4}+\frac{1}{128}c_i^{14}z^{-10/4}+\cdots,$$
and $c_1=1$, $c_2=I$, $c_3=-1$,  $c_4=-I$. Note that the first term of
these series is $z^{1/2}$ or $-z^{1/2}$.  The curve $\overline{{\cal C}}$ also has one infinity branch defined by  $\overline{B}=\bigcup_{i=1}^2 \overline{L}_i$, where
$\overline{L}_i=\{(z,\overline{r}_i(z))\in {\Bbb C}^2: \,z\in {\Bbb
C},\,|z|>\overline{M}\}$, $\overline{r}_i(z)=d_i z^{1/2},$ and
$d_1=1$, $d_2=-1$. We get that $B$ and $\overline{B}$, are convergent
since $L_1$ and $\overline{L}_1$ converge. In fact, $L_1$ and $L_3$
converge with $\overline{L}_1$ and, on the other hand, $L_2$ and
$L_4$ converge with $\overline{L}_2$ (see Lemma
\ref{L-DistVertical}).
\end{remark}

\para

Two convergent infinity branches may be contained in the same curve
or they may belong to different curves. In this second case we will
say that those curves {\it approach each other}. In order to define
this concept in a more formal way, we first introduce the following
distance:

\begin{definition}\label{D-distance}
Given an algebraic plane curve ${\cal C}$ over $\Bbb C$ and a point $p\in {\Bbb
C}^2$, we define {\em the distance from $p$ to ${\cal C}$} as
$ d(p,{\cal C})=\min\{d(p,q):q\in{\cal C}\}.$
\end{definition}

\begin{remark}
Observe that this minimum exists because ${\cal C}$ is a closed set.
\end{remark}

\begin{definition}\label{D-distance1}
Let ${\cal C} $ be an algebraic plane curve over ${\Bbb C}$ with an
infinity branch $B$. We say that a  curve ${\overline{{\cal C}}}$
{\em approaches} ${\cal C}$ at its infinity branch $B$ if there
exists one leaf $L=\{(z,r(z))\in {\Bbb C}^2:\,z\in {\Bbb
C},\,|z|>M\}\subset B$ such that
$\lim_{z\rightarrow\infty}d((z,r(z)),\overline{\cal C})=0.$
\end{definition}

We will show that this condition is satisfied for one leaf of
$B$ if and only if it is satisfied for every leaf of $B$. It will be
derived as a consequence of the following theorem.

\para

  \begin{theorem}\label{T-curvas-aprox} Let ${\cal C}$
be a plane algebraic curve over $\Bbb C$ with an infinity branch
$B$. A plane algebraic curve ${\overline{{\cal C}}}$ approaches
${\cal C}$ at $B$ if and only if ${\overline{{\cal C}}}$ has an
infinity branch, $\overline{B}$, such that $B$ and $\overline{B}$
are convergent.
\end{theorem}

\noindent\textbf{Proof:} Suppose that ${\overline{{\cal C}}}$
approaches ${\cal C}$ at $B$. Then, there exists a leaf
$L=\{(z,r(z))\in {\Bbb C}^2:\,z\in {\Bbb C},\,|z|>M\}\subset B$ such
that
$\lim_{z\rightarrow\infty}d((z,r(z)),\overline{\cal C})=0.$
In addition, let $P=(1:m:0)$ be the infinity point associated to
$B$, and let $\{z_n\}_{n\in \mathbb{N}}$ be a sequence in ${\Bbb C}$
such that $\lim_{n\rightarrow\infty} z_n= \infty$.  We have that
$\lim_{n\rightarrow\infty}d((z_n,r(z_n)),\overline{\cal C})=0$ which
implies that
$$\lim_{n\rightarrow\infty}d((z_n,r(z_n)),(p_n,q_n))=0,$$ where, for
each $z_n$ such that $|z_n|>M$, $(p_n,q_n)$ is the point of ${\overline{{\cal C}}}$ closest to the point $(z_n,r(z_n))$ (this point exists because
of Definition \ref{D-distance}). Note that the above equality
implies that
$$
\lim_{n\rightarrow\infty}|p_n-z_n|^2+|q_n-r(z_n)|^2=0,
$$
and hence we have that:

\begin{itemize}
\item
$\displaystyle\lim_{n\rightarrow\infty}(p_n-z_n)=0$. Then, $\lim_{n\rightarrow\infty} p_n/z_n=1$ which implies that $\lim_{n\rightarrow\infty} p_n=\infty$. Hence,  $\lim_{n\rightarrow\infty}1/p_n=0.$
\item $\displaystyle\lim_{n\rightarrow\infty}(q_n-r(z_n))=0$. Then, $\lim_{n\rightarrow\infty}(q_n/z_n-r(z_n)/z_n)=0$ which implies that
$\lim_{n\rightarrow\infty} q_n/z_n =\lim_{n\rightarrow\infty} r(z_n)/z_n =m.$
\end{itemize}
Therefore,
$$\lim_{n\rightarrow\infty} q_n/p_n =\lim_{n\rightarrow\infty} \frac{q_n/z_n}{p_n/z_n} =m.$$
Now, taking into account Lemma \ref{L-BolaMaxima} and  that
$\lim_{n\rightarrow\infty} p_n=\infty$, we get that there exits
$n_0\in {\Bbb N}$ such that for $n\geq n_0$, the points $(p_n,q_n)$
are in some infinity branch of ${\overline{{\cal C}}}$. Moreover,
since any curve has a finite number of infinity branches and a
finite number of leaves, we can find a subsequence
$\{z_{n_l}\}_{l\in {\Bbb N}}$ and $l_0\in {\Bbb N}$ such that for
$l\geq l_0$,   the points $(p_{n_l},q_{n_l})$ are all in a same leaf
$\overline{L}=\{(z,\overline{r}(z))\in {\Bbb C}^2:\,z\in {\Bbb
C},\,|z|>\overline{M}\}$, belonging to some branch
$\overline{B}\subset{\overline{{\cal C}}}$.
\para

\noindent
Under these conditions, we deduce that for $l\geq l_0$, $q_{n_l}=\overline{r}(p_{n_l})$, and then
$$\displaystyle\lim_{l\rightarrow\infty}\overline{r}(p_{n_l})/p_{n_l}=\lim_{l\rightarrow\infty}q_{n_l}/p_{n_l}=m.$$
Since the limit  $\lim_{z\rightarrow\infty}\frac{\overline{r}(z)}{z}=\lim_{z\rightarrow\infty}\overline{\varphi}(z^{-1})$ exists, we get
$\lim_{z\rightarrow\infty}\frac{\overline{r}(z)}{z}=m$.

\para

\noindent In addition, note that
$$|r(z_{n_l})-\overline{r}(z_{n_l})|=d((z_{n_l},r(z_{n_l})),(z_{n_l},\overline{r}(z_{n_l})))\leq$$
$$d((z_{n_l},r(z_{n_l})),(p_{n_l},\overline{r}(p_{n_l})))
+d((p_{n_l},\overline{r}(p_{n_l})),(z_{n_l},\overline{r}(z_{n_l})))\qquad \mbox{(I)}$$
and
$$d((z_{n_l},r(z_{n_l})),(p_{n_l},\overline{r}(p_{n_l})))=d((z_{n_l},r(z_{n_l})),(p_{n_l},q_{n_l}))_{\overrightarrow{l\rightarrow
\infty}}0.$$ Now, let us prove that
$d((p_{n_l},\overline{r}(p_{n_l})),(z_{n_l},\overline{r}(z_{n_l})))_{\overrightarrow{l\rightarrow
\infty}}0.$
For this purpose, we show that
$\lim_{l\rightarrow\infty} (\overline{r}(p_{n_l})-\overline{r}(z_{n_l}))=0$. Indeed: let $$\overline{r}(z)=mz+b_1z^{\frac{s-s_1}{s}}+b_2z^{\frac{s-s_2}{s}}+\cdots,\,\,\, s, s_i\in {\Bbb N},\, 0<s_1<s_2<\cdots.$$
Thus,
$$\overline{r}\,'(z)={m} +\frac{s-s_1}{s}b_1z^{\frac{-s_1}{s}}+\frac{s-s_2}{s}b_2z^{\frac{-s_2}{s}}+\cdots\,_{\overrightarrow{z\rightarrow
\infty}} \,{m}.$$ Therefore, there exist $K>0$ and $\delta>0$ such that
$|\overline{r}\,'(z)|\leq K$, for $|z|>\delta$. Applying the Mean Value
Theorem (see \cite{Ahlfors}), we have that
$$\Re\left(\frac{\overline{r}(p_{n_l})-\overline{r}(z_{n_l})}{p_{n_l}-z_{n_l}}\right)=\Re (\overline{r}\,'(c_1)),\qquad
\Im\left(\frac{\overline{r}(p_{n_l})-\overline{r}(z_{n_l})}{p_{n_l}-z_{n_l}}\right)=\Im (\overline{r}\,'(c_2)),$$
where $\Re(q)$ and $\Im(q)$ denote the real part and the imaginary part of $q(z)\in {\Bbb C}(z)$,   respectively,  and $c_1, c_2\in ]p_{n_l}, z_{n_l}[$, where $]p_{n_l}, z_{n_l}[:= \{z \in {\Bbb C}:\, z =p_{n_l}+ (p_{n_l} - z_{n_l})t,\, t\in (0,1)\}.$ Thus,
$$|\overline{r}(p_{n_l})-\overline{r}(z_{n_l})|^2=(\Re (\overline{r}\,'(c_1))^2+\Im (\overline{r}\,'(c_2))^2)|p_{n_l}-z_{n_l}|^2.$$
Now,
since $\lim_{l\rightarrow\infty}p_{n_l}=\lim_{l\rightarrow\infty}z_{n_l}=\infty$, and $\lim_{l\rightarrow\infty}p_{n_l}-z_{n_l}=0$, we deduce that   given
  $\varepsilon>0$, there exists  $l_1\in {\Bbb N}$ such that, for $l\geq l_1$,   $$|p_{n_l}|>\delta+\varepsilon,\quad |z_{n_l}|>\delta+\varepsilon, \quad \mbox{and}\quad
|p_{n_l}-z_{n_l}|<\varepsilon.$$ Then,
$|c_j|>\delta$ and   $|\overline{r}\,'(c_j)|\leq K$ for $j=1,2$, which implies that,
for $l\geq l_1$,
$$|\overline{r}(p_{n_l})-\overline{r}(z_{n_l})|\leq \sqrt{2}K|p_{n_l}-z_{n_l}|_{\overrightarrow{l\rightarrow
\infty}}0 $$
(note that $\Re (\overline{r}\,'(c_1))\leq |\overline{r}\,'(c_1)|\leq K$, and $\Im (\overline{r}\,'(c_2))\leq |\overline{r}\,'(c_2)|\leq K$).
Therefore, $\lim_{l\rightarrow\infty} (\overline{r}(p_{n_l})-\overline{r}(z_{n_l}))=0$, which implies that there exists a sequence
$\{z_{n_l}\}_{l\in {\Bbb N}}$ with
 $\lim_{l\rightarrow\infty} z_{n_l}=\infty$, such that
$$\lim_{l\rightarrow\infty}(r(z_{n_l})-\overline{r}(z_{n_l}))=0 $$
(see inequality (I)). Then,  the terms with positive exponent of the
series $r(z)$ and $\overline{r}(z)$ are the same (see the proof of
Lemma \ref{L-DistVertical}). Hence, we conclude that   (see Lemma \ref{L-DistVertical})
$$\lim_{z\rightarrow\infty}(r(z)-\overline{r}(z))=0$$
and thus $B$ and $\overline{B}$ are convergent (see Definition
\ref{D-distance0}).

\para

Reciprocally, let us assume that $B$ and $\overline{B}$ are convergent.
Then, by definition, there exist two leaves $L=\{(z,r(z))\in {\Bbb
C}^2:\,z\in {\Bbb C},\,|z|>M\}\subset B$ and
$\overline{L}=\{(z,\overline{r}(z))\in {\Bbb C}^2:\,z\in {\Bbb
C},\,|z|>\overline{M}\}\subset \overline{B}$ such that
$\displaystyle\lim_{z\rightarrow\infty}(r(z)-\overline{r}(z))=0$. Therefore,
$$\lim_{z\rightarrow\infty}d((z,r(z)),\overline{\cal C})\leq
\lim_{z\rightarrow\infty}d((z,r(z)),(z,\overline{r}(z)))=\lim_{z\rightarrow\infty}(r(z)-\overline{r}(z))=0.\qquad \mbox{\hfill $\Box$}$$

\para
\begin{remark}
\begin{enumerate}
\item From Theorem \ref{T-curvas-aprox}, we get that ``{\em proximity}'' is a
symmetric relation; i.e., ${\overline{{\cal C}}}$ approaches ${\cal C}$ at
some infinity branch $B$ iff ${\cal C}$ approaches
${\overline{{\cal C}}}$ at some infinity branch $\overline{B}$. In the following, we
say that ${\cal C}$ and ${\overline{{\cal C}}}$ approach each other
or that they are {\em approaching curves}.
\item Theorem \ref{T-curvas-aprox} and Remark
\ref{R-common-inf-point} imply that two approaching curves have
 a common infinity point.
\item From Theorem \ref{T-curvas-aprox} and Proposition
\ref{P-convergent-branches}, we get that ${\overline{{\cal C}}}$ {\em
approaches} ${\cal C}$ at an infinity branch $B$ if and only if for every leaf $L=\{(z,r(z))\in {\Bbb C}^2:\,z\in {\Bbb
C},\,|z|>M\}\subset B$, it holds that
$\displaystyle\lim_{z\rightarrow\infty}d((z,r(z)),\overline{\cal
C})=0$.
\end{enumerate}
\end{remark}

\begin{corollary}\label{C-approaching-curves}
Let $\cal C$ be an algebraic plane curve with an infinity branch
$B$. Let ${\overline{{\cal C}}}_1$ and ${\overline{{\cal C}}}_2$ be
two different curves that approach $\cal C$ at $B$. Then
${\overline{{\cal C}}}_1$ and ${\overline{{\cal C}}}_2$ approach
each other.
\end{corollary}

\noindent\textbf{Proof:}  From Theorem
\ref{T-curvas-aprox}, there
 exist two infinity branches
$B_1\subset {\overline{{\cal C}}}_1$ and $B_2\subset
{\overline{{\cal C}}}_2$, convergent with $B$. Thus, for
each leaf $L=\{(z,r(z))\in {\Bbb C}^2:\,z\in {\Bbb
C},\,|z|>M\}\subset B$, there exist two leaves $L_1=\{(z,r_1(z))\in
{\Bbb C}^2:\,z\in {\Bbb C},\,|z|>M_1\}\subset B_1$ and
$L_2=\{(z,r_2(z))\in {\Bbb C}^2:\,z\in {\Bbb C},\,|z|>M_2\}\subset
B_2$  such that $\lim_{z\rightarrow\infty}(r(z)-r_1(z))=0$ and
$\lim_{z\rightarrow\infty}(r(z)-r_2(z))=0$. Then
$$|r_1(z)-r_2(z)|\leq |r_1(z)-r(z)|+|r(z)-r_2(z)|_{\overrightarrow{z\rightarrow
\infty}}0.$$  Therefore, ${\overline{{\cal C}}}_1$ and ${\overline{{\cal C}}}_2$ approach each other.\hfill $\Box$

\para

\para

\noindent
In the following, we illustrate the above results with an example.

\para

\begin{example}
Let $\cal C$ and $\overline{{\cal C}}$ be two plane curves defined implicitly by the polynomials
$$f(x,y)=2y^3x-y^4+2y^2x-y^3-2x^3+x^2y+3\in {\Bbb R}[x,y],\qquad \mbox{and}$$
$$\overline{f}(x,y)=y^3x-y^4+y^2x-y^3-x^3+x^2y+2\in {\Bbb R}[x,y],$$
respectively. Let us prove that   $\cal C$ and $\overline{{\cal C}}$ approach each other (see Figure
\ref{Ej-approaching curves}) at the infinity branch  associated to
the infinity point $P=(1:0:0)$ (note that both curves have $P$ as an
infinity point). Reasoning as in Example \ref{Ej-ramas infinitas},
we get that the infinity branch of $\cal C$ associated to $P$ is
given by $B=L_1\cup L_2\cup L_3$, where $L_i=\{(z,r_i(z))\in {\Bbb
C}^2:\,z\in {\Bbb C},\,|z|>M\}$,
$$r_i(z)=c_i^2z^{2/3}-1/3+1/9c_i^2z^{-2/3}-2/81c_i^4z^{-4/3}-1/2c_i^7z^{-7/3}+\cdots$$
and $c_i,\,i=1,2,3$ are the complex roots of $x^3=1$. On the other
hand, the infinity branch of $\overline{{\cal C}}$ associated to $P$
is given by $\overline{B}=\overline{L}_1\cup \overline{L}_2\cup
\overline{L}_3$, where $\overline{L}_i=\{(z,\overline{r}_i(z))\in
{\Bbb C}^2:\,z\in {\Bbb C},\,|z|>M\}$,
$$\overline{r}_i(z)=c_i^2z^{2/3}-1/3+1/9c_i^2z^{-2/3}-2/81c_i^4z^{-4/3}-2/3c_i^7z^{-7/3}+\cdots$$
and $c_i,\,i=1,2,3$ are the complex roots of $x^3=1$ (to compute $r_i$ and $\overline{r}_i$, we use  the {\it
algcurves} package included in {\it Maple}). From Lemma
\ref{L-DistVertical}, we conclude that both branches converge, since
the terms with non negative exponent  in both series, $r_i$ and
$\overline{r}_i$, are the same.
\vspace*{-1cm}
\begin{figure}[h]
$$
\begin{array}{lcr}
\psfig{figure=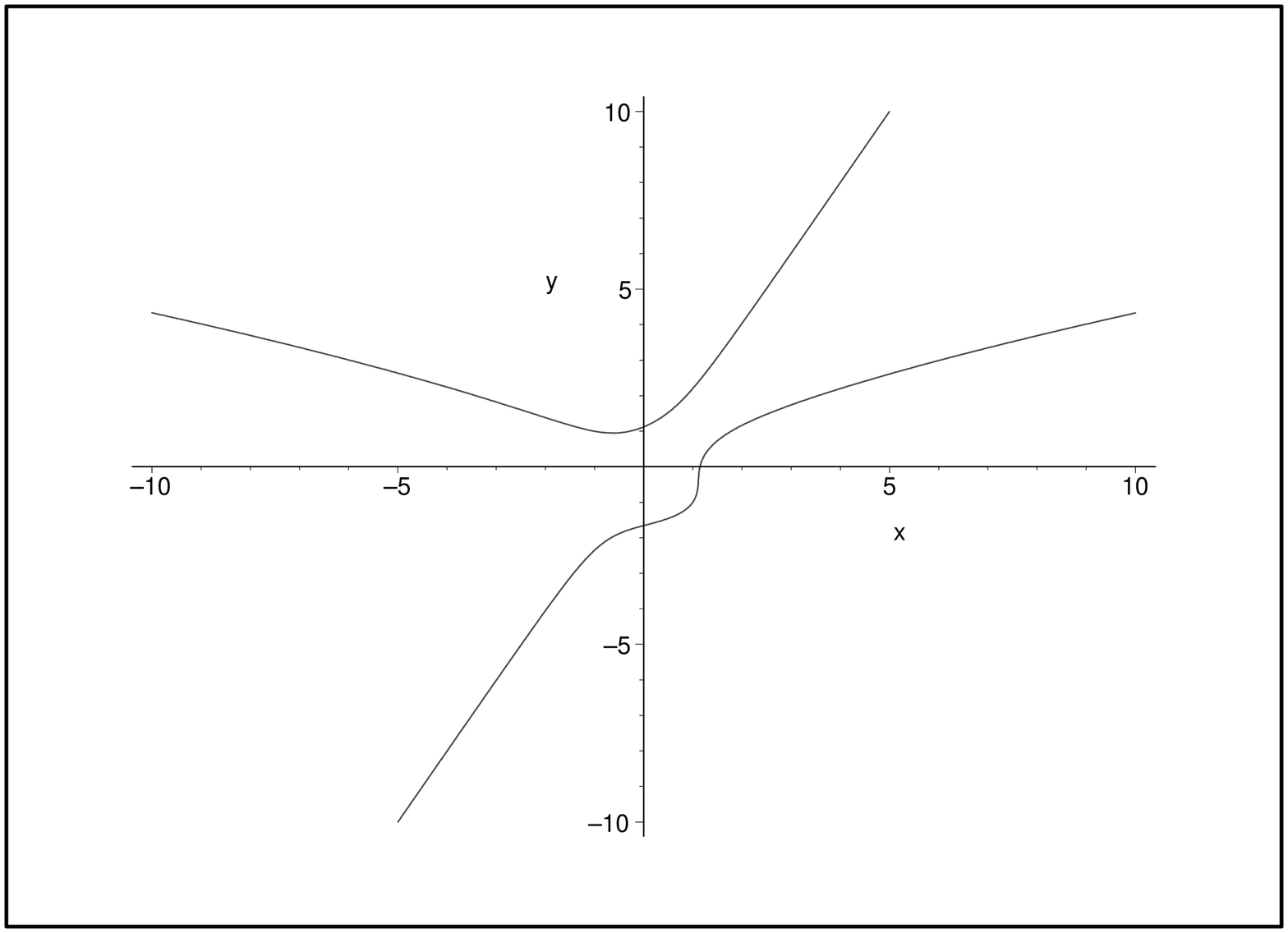,width=4.3cm,height=4.3cm,angle=270} &
\psfig{figure=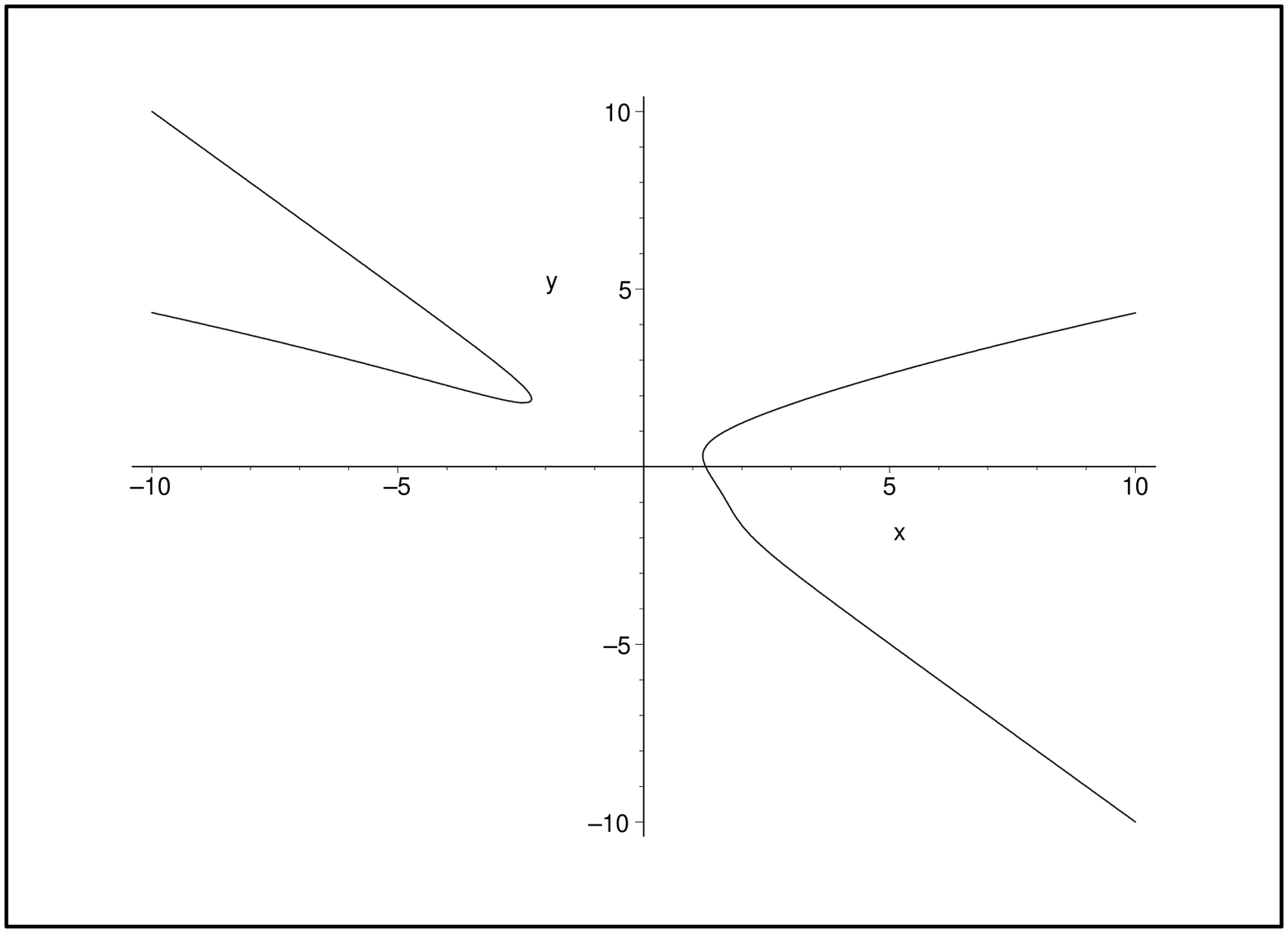,width=4.3cm,height=4.3cm,angle=270} &
\psfig{figure=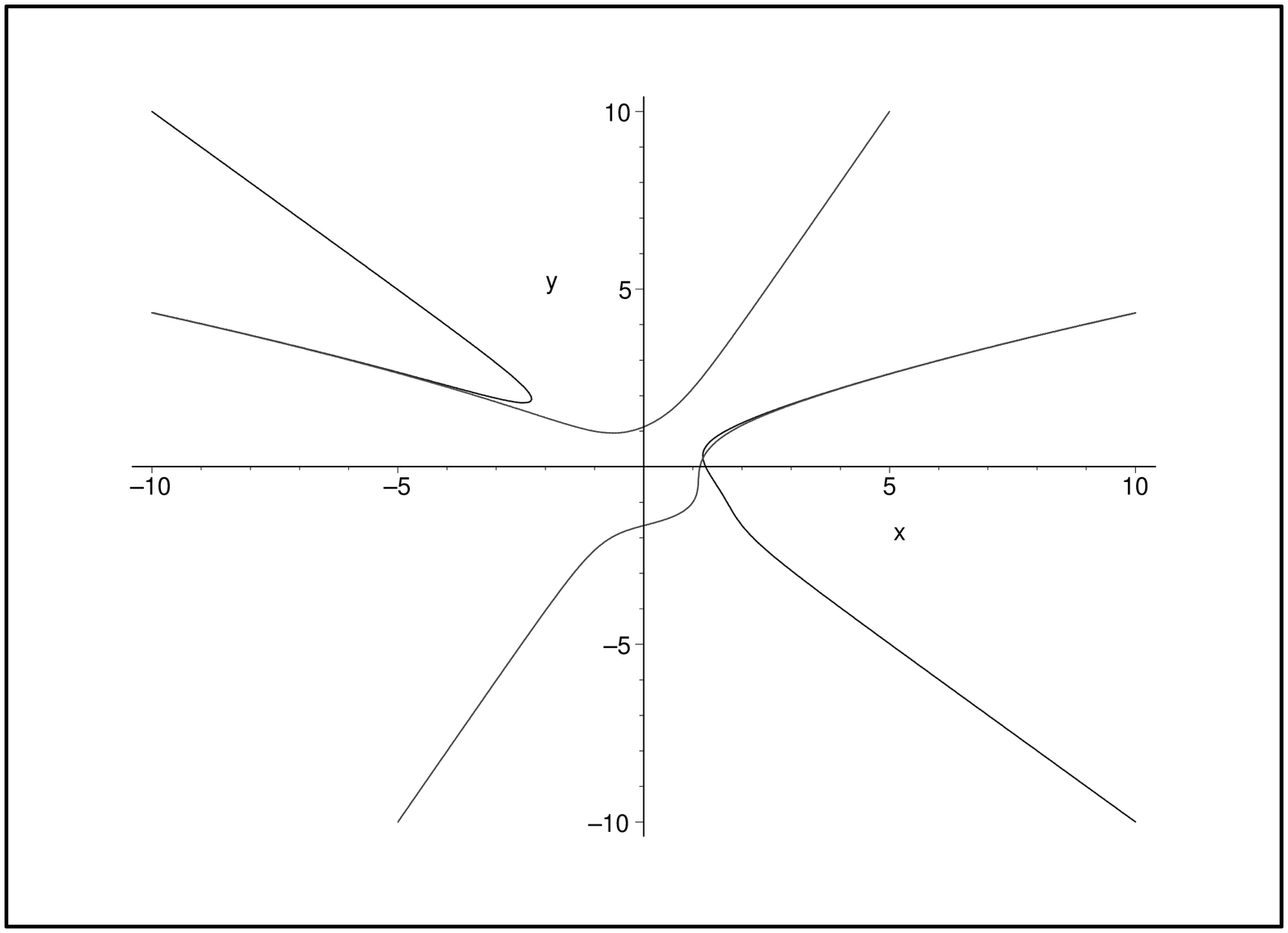,width=4.3cm,height=4.3cm,angle=270}
\end{array}
\vspace*{-0.5cm}
$$ \caption{$\cal C$ (left),  $\overline{{\cal C}}$ (center), and both approaching curves (right)}\label{Ej-approaching curves}
\end{figure}
\end{example}

\section{Asymptotic Behavior}

Using the results presented in the previous sections, in the following we analyze the behavior of two curves at the infinity (the {\it asymptotic behavior}). More precisely, in this
section we present an algorithm that provides a method to compare the
behavior of two algebraic plane curves  as they tend to infinity. In addition, we prove that if two plane algebraic curves have
the same {\it asymptotic behavior},   the Hausdorff distance between
them is finite.

\para

\noindent
To start with, we first introduce the following definition.

\begin{definition} We say that two algebraic plane
curves, ${\cal C}$ and ${\overline{{\cal C}}}$, have the same asymptotic behavior
if every infinity branch of ${\cal C}$ converges to another branch
of ${\overline{{\cal C}}}$, and reciprocally.
%\begin{enumerate}
%\item[a)] for each infinity branch $B$ of ${\cal C}$, ${\overline{{\cal C}}}$
%has an infinity branch convergent to $B$.
%\item[b)] for each infinity branch $\overline{B}$ of ${\overline{{\cal C}}}$, ${\cal C}$
%has an infinity branch convergent to $\overline{B}$.
%\end{enumerate}
\end{definition}

\begin{remark}\label{R-same-asymp-behav} From Theorem \ref{T-curvas-aprox}, we deduce that ${\cal C}$
and ${\overline{{\cal C}}}$ have the same asymptotic behavior if and only if ${\cal C}$
approaches ${\overline{{\cal C}}}$ at all its infinity branches, and reciprocally.
\end{remark}

\para

\noindent Now, we recall the notion of  {\em Hausdorff distance}.

\begin{definition}\label{D-Hausdorff} Given a metric space $(E,d)$ and two
subsets $A, B\subset E\setminus \{\emptyset\}$, the {\em Hausdorff distance} between them is defined
as:
$$d_H(A,B)=\max\{\sup_{x\in A}\inf_{y\in B}d(x,y),\sup_{y\in B}\inf_{x\in A}d(x,y)\}.$$
If $E=\mathbb{C}^2$ and $d$ is the Euclidean distance, the Hausdorff distance between two curves ${\cal C}$ and
${\overline{{\cal C}}}$ can be expressed as:
$$d_H({\cal C},{\overline{{\cal C}}})=\max\{\sup_{p\in {\cal C}}d(p,{\overline{{\cal C}}}),\sup_{\overline{p}\in {\overline{{\cal C}}}}d(\overline{p},{\cal C})\}.$$
%where $d(p,{\overline{{\cal C}}})$ is the point-curve distance given in Definition \ref{D-distance}.
\end{definition}

\begin{proposition}\label{P-Hausdorff}
Let ${\cal C}$ and ${\overline{{\cal C}}}$ be two algebraic plane curves having
the same asymptotic behavior. Then, the Hausdorff distance between
them is finite.
\end{proposition}

\noindent\textbf{Proof:} Let $r$ be the number of infinity branches
of ${\cal C}$. Then,
${\cal C}=B_1\cup \cdots \cup B_r \cup \widehat{B},$
where $\widehat{B}$ is the set of points of ${\cal C}$ that do not
belong to any infinity branch. Thus,
$$\sup_{p\in {\cal C}}d(p,{\overline{{\cal C}}})=\max\{\sup_{p\in {B_1}}d(p,{\overline{{\cal C}}}),...,\sup_{p\in {B_r}}
d(p,{\overline{{\cal C}}}),\sup_{p\in
\widehat{B}}d(p,{\overline{{\cal C}}})\}.$$ For each
$i=1,...,r$, let $B_i=\bigcup_{j=1}^{N_i}L_{i,j}$, where
$L_{i,j}=\{(z,r_{i,j}(z))\in {\Bbb C}^2:\,z\in {\Bbb C},\,|z|>M_i\}$, and
$N_i=\nu(B_i)$. Then,
$$\sup_{p\in {B_i}}d(p,{\overline{{\cal C}}})=\max_{j=1,\ldots,N_i}\left\{\sup_{|z|>M_i}d((z,r_{i,j}(z)),{\overline{{\cal C}}})\right\}.$$
Moreover, from Remark \ref{R-same-asymp-behav},
${\overline{{\cal C}}}$ approaches ${\cal C}$ at $B_i$, so
$\lim_{z\rightarrow\infty}d((z,r_{i,j}(z)),\overline{\cal C})=0$ for
every $j=1,\ldots,N_i$. Hence, given   $\varepsilon>0$ there exists
$\delta>0$ such that $d((z,r_{i,j}(z)),\overline{\cal
C})<\varepsilon$, for $|z|>\delta$.
%\textcolor{red}{(esto habria que verlo despacio, nuestra definicion
%de limite no dice exactamente esto)}.
Then, since $r_{i,j}$ is a continuous function, and $\{z\in {\Bbb
C}:\,M_i\leq|z|\leq\delta\}$ is a compact set, we deduce that
$$\sup_{p\in {B_i}}d(p,{\overline{{\cal C}}})\leq \max_{j=1,\ldots,N_i}\max\left\{\sup_{M_i\leq|z|\leq\delta}d((z,r_{i,j}(z)),
{\overline{{\cal C}}}),\varepsilon\right\}<\infty.$$

\para

Now, let  $p=(a,b)\in \widehat{B}$. From Lemma
\ref{L-BolaMaxima} and Remark \ref{R-BolaMaxima}, we have that there exists  $K\in \mathbb{R}^+$
such that $|a|,|b|\leq K$. Thus, $d(p,{\cal O})\leq K$, where ${\cal O}$ is the
origin and, $$d(p,{\overline{{\cal C}}})\leq d(p,{\cal O})+d({\cal O},{\overline{{\cal C}}})\leq K+d({\cal O},{\overline{{\cal C}}}).$$ Note that $K<\infty$, and $d({\cal O},{\overline{{\cal C}}})<\infty,$  which implies that $\sup_{p\in \widehat{B}}d(p,{\overline{{\cal C}}})<\infty$.

\para

Therefore, we conclude that $\sup_{p\in {\cal C}}d(p,{\overline{{\cal C}}})<\infty$. Reasoning similarly, we deduce that $\sup_{\overline{p}\in
{\overline{{\cal C}}}}d(\overline{p},{\cal C})<\infty$, which implies that $d_H({\cal C},\overline{\cal C})<\infty$.
 \hfill $\Box$

\para

\para

The following algorithm allow us to compare the asymptotic
behavior of two   curves $\cal C$ and $\overline{{\cal C}}$.  We assume that we have  prepared   $\cal C$ and $\overline{{\cal C}}$ such that  by means  of a suitable linear change of coordinates (the same change applied to both curves),  $(0:1:0)$  is not a point of infinity of ${\cal C}^*$ and $\overline{{\cal C}}^*$.
%neither ${\cal C}^*$  nor $\overline{{\cal C}}^*$ has $(0:1:0)$  as a point of infinity.

\para

\noindent
\begin{center}
\fbox{\hspace*{2 mm}\parbox{13.2cm}{ \vspace*{2 mm} {\bf Algorithm
{\sf Asymptotic Behavior.}}
\vspace*{0.2cm}

\noindent
Given two implicit algebraic plane curves $\cal C$ and $\overline{{\cal C}}$, the algorithm
decides whether $\cal C$ and $\overline{{\cal C}}$ have the same asymptotic behavior.
\begin{enumerate}
\item[1.] Compute the infinity points of $\cal C$ and $\overline{{\cal C}}$. If they are
not the same,  {\sc Return} {\it the curves do not have the same
asymptotic behavior} (see Remark \ref{R-common-inf-point}). Otherwise, let $P_1,\ldots,P_n$ be these infinity points.
\item[2.]
For each $P_k:=(1:m_k:0)$,\,$k=1,\ldots, n$ do:
\begin{enumerate}
\item[2.1.] Compute the infinity branches of $\cal C$ associated
to $P_k$. Let $B_1,...,B_{n_k}$ be these branches. For
each $i=1,\ldots,n_k$, let $L_i=\{(z,r_i(z))\in {\Bbb C}^2:\,z\in
{\Bbb C},\,|z|>M_i\}$ be any leaf of $B_i$.% randomly chosen
\item[2.2.] Compute the infinity branches  of $\overline{{\cal C}}$ associated
to $P_k$. Let
$\overline{B}_1,...,\overline{B}_{l_k}$ be these branches. For
each $j=1,\ldots,l_k$, let
$\overline{L}_j=\{(z,\overline{r}_j(z))\in {\Bbb C}^2:\,z\in {\Bbb
C},\,|z|>M_j\}$ be any leaf of $\overline{B}_j$.% randomly chosen
\item[2.3.] For each $B_i\subset {\cal C}$, find
  $\overline{B}_j\subset{\overline{{\cal C}}}$ such that the
terms with non negative exponent in $r_i(z)$ and $\overline{r}_j(z)$
are the same up to conjugation. If there isn't such a branch,  {\sc
Return} {\it the curves do not have the same asymptotic behavior}
(see Lemma \ref{L-DistVertical}).
\item[2.4.] For each  $\overline{B}_j\subset{\overline{{\cal C}}}$, find
   $B_i\subset {\cal C}$ such that the
terms with non negative exponent in $r_i(z)$ and $\overline{r}_j(z)$
are the same up to conjugation. If there isn't such a branch,  {\sc
Return} {\it the curves do not have the same asymptotic behavior}
(see Lemma \ref{L-DistVertical}).
\end{enumerate}
\item[3.] {\sc Return} {\it the curves $\cal C$ and $\overline{{\cal C}}$ have the same asymptotic behavior}.
\end{enumerate}}\hspace{2 mm}}
\end{center}

\para

In the following, we illustrate  the performance of algorithm
{\sf Asymptotic Behavior} with an example.

\para

\begin{example}
Let $\cal C$, and $\overline{{\cal C}}$ be two plane curves defined implicitly by the polynomials
$$f(x,y)=2y^3x-y^4+2y^2x-y^3-2x^3+x^2y+3,\qquad \mbox{and}$$
$$\overline{f}(x,y)=2y^3x-y^4+2y^2x-y^3-2x^3+x^2y-3x^2-xy+2x-3y+1,$$
respectively. We apply  the algorithm {\sf Asymptotic
Behavior} to decide whether $\cal C$ and $\overline{{\cal C}}$  have the same asymptotic
behavior:

\begin{itemize}
\item[] \mbox{\sf {Step 1}:}  Compute
the infinity points of  $\cal C$ and $\overline{{\cal C}}$. We obtain that $\cal C$  and $\overline{{\cal C}}$  have the same infinity
points: $P_1=(1:0:0)$ and $P_2=(1:2:0)$.
\para

We start by analyzing  the infinity branches associated to $P_1$:

\item[]  \mbox{\sf {Step 2.1}:} Reasoning as in Example \ref{Ej-ramas infinitas}, we get that the only infinity branch associated to $P_1$
in $\cal C$ is given by $B_1=L_{1,1}\cup L_{1,2}\cup L_{1,3}$ where
$L_{1,i}=\{(z,r_{1,i}(z))\in {\Bbb C}^2:\,z\in {\Bbb C},\,|z|>M_1\}$,
$i=1,2,3$, and
$$r_{1,i}(z)=z^{2/3}-1/3+1/9z^{-2/3}-2/81z^{-4/3}+\cdots,$$ up to
conjugation.
\item[]  \mbox{\sf {Step 2.2}:} We also have that there exists only one infinity branch associated to $P_1$ in $\overline{{\cal C}}$. It  is given by
$\overline{B}_1=\overline{L}_{1,1}\cup \overline{L}_{1,2}\cup
\overline{L}_{1,3}$ where
$\overline{L}_{1,i}=\{(z,\overline{r}_{1,i}(z))\in {\Bbb C}^2:\,z\in
{\Bbb C},\,|z|>\overline{M}_1\}$, $i=1,2,3$, and
$$\overline{r}_{1,i}(z)=z^{2/3}-1/3+1/2z^{-1/3}+19/36z^{-2/3}+\cdots,$$
up to conjugation.
\item[]  \mbox{\sf {Step 2.3}} and   \mbox{\sf {Step 2.4}:}
$r_{1,1}(z)$ and $\overline{r}_{1,1}(z)$ have the same terms with non negative exponent. Thus, ${B}_1$ and $\overline{B}_1$
  converge.\\

Now we analyze the infinity branches associated to $P_2$:

\item[]  \mbox{\sf  {Step 2.1}:} Reasoning as in Example \ref{Ej-ramas infinitas}, we get that the only infinity branch associated
to $P_2$ in $\cal C$   is given by $B_2=L_2=\{(z,r_2(z))\in {\Bbb
C}^2:\,z\in {\Bbb C},\,|z|>M_2\}$, where

$$r_2(z)=2z+3/8z^{-3}-9/64z^{-4}+27/512z^{-5}+\cdots.$$

\item[]  \mbox{\sf  {Step 2.2}:} The only infinity branch associated to $P_2$ in $\overline{{\cal C}}$  is given by
$\overline{B}_2=\overline{L}_2=\{(z,\overline{r}_2(z))\in {\Bbb
C}^2:\,z\in {\Bbb C},\,|z|>\overline{M}_2\}$, where

$$\overline{r}_2(z)=2z-5/8z^{-1}-17/64z^{-2}-145/512z^{-3}+\cdots.$$

\item[]  \mbox{\sf  {Step 2.3}} and   \mbox{\sf {Step 2.4}:}
$r_2(z)$ and $\overline{r}_2(z)$ have the same terms with non negative exponent. Thus, ${B}_2$ and $\overline{B}_2$  converge.
\end{itemize}

Since every infinity branch of $\cal C$ converges to another branch
of $\overline{{\cal C}}$, and reciprocally, the algorithm returns that $\cal C$ and $\overline{{\cal C}}$ have the
same asymptotic behavior (see Figure \ref{Ej-comportamiento
asintotico}).
\vspace*{-1cm}
\begin{figure}[h]
$$
\begin{array}{lcr}
\psfig{figure=Ejemplo2a.eps,width=4.3cm,height=4.3cm,angle=270} &
\psfig{figure=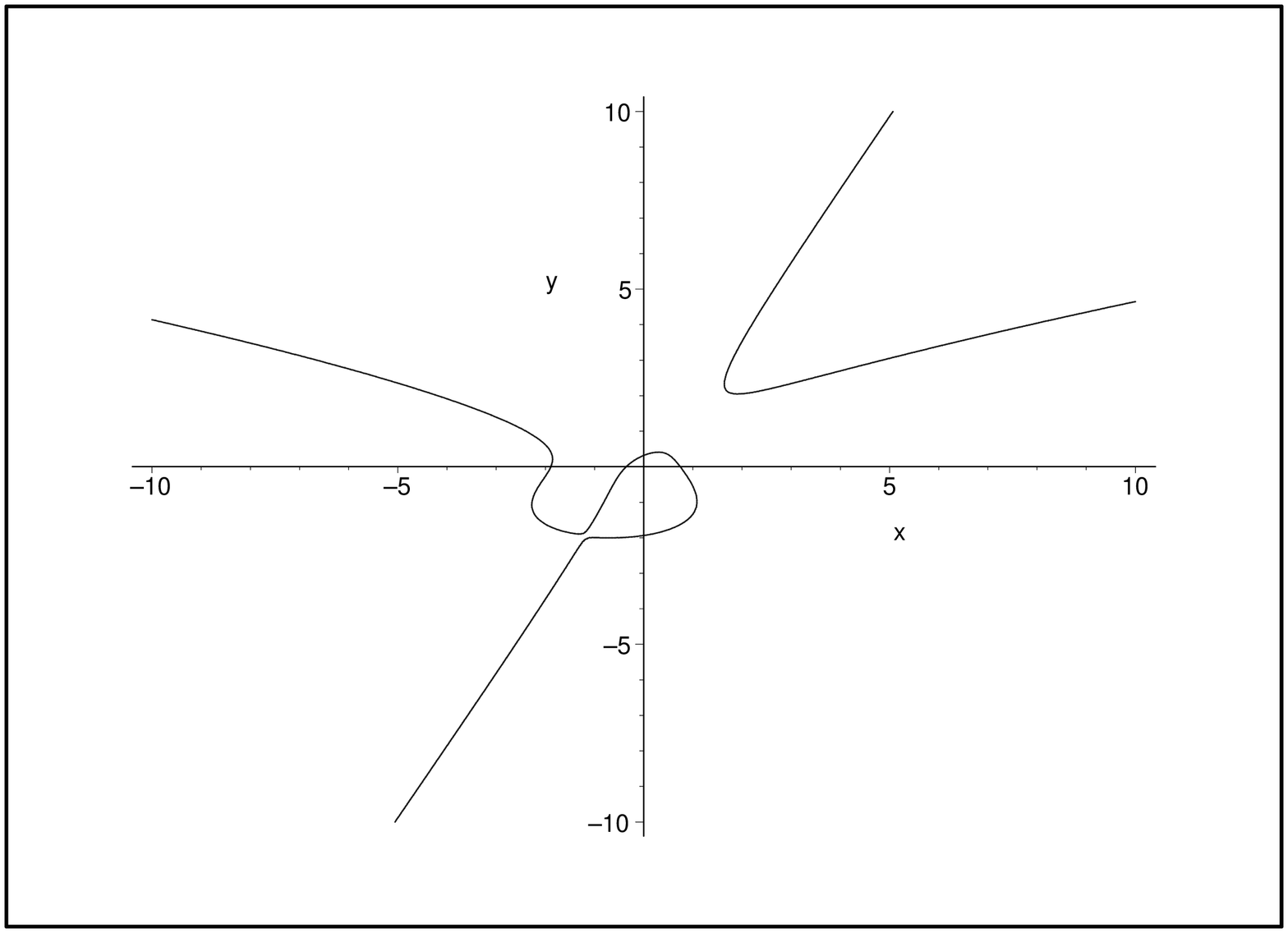,width=4.3cm,height=4.3cm,angle=270} &
 \psfig{figure=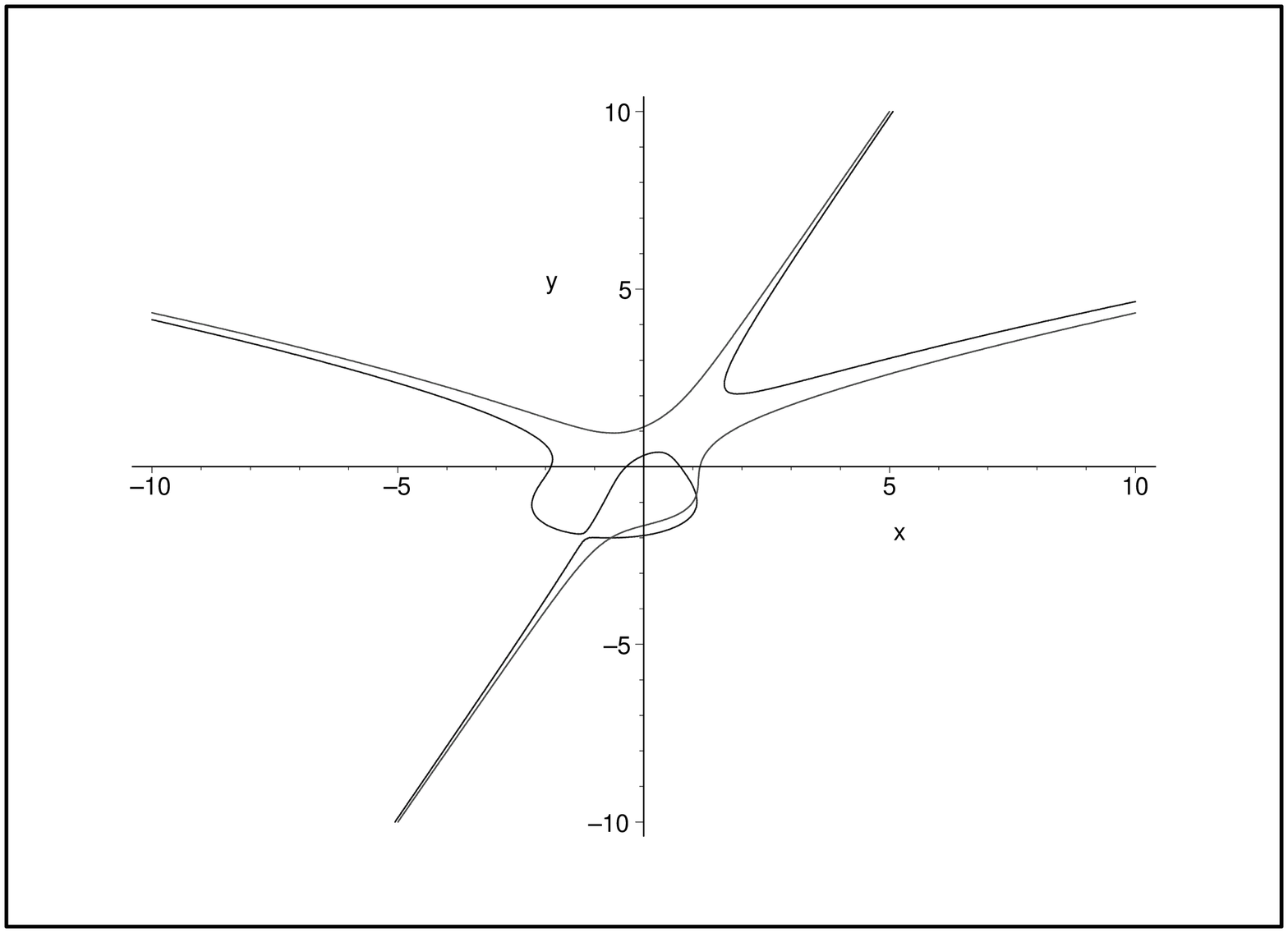,width=4.3cm,height=4.3cm,angle=270}\end{array}
 \vspace*{-0.5cm}
$$ \caption{$\cal C$ (left),  $\overline{{\cal C}}$ (center), and the asymptotic behavior of $\cal C$ and $\overline{{\cal C}}$ (right)}\label{Ej-comportamiento asintotico}
\end{figure}
\end{example}


\begin{thebibliography}{00}


 \bibitem{Ahlfors} Ahlfors, L.V.  (1979). {\it Complex Analysis}.  McGraw-Hill, Third Edition.


\bibitem{paper2} Blasco, A., P\'erez-D\'{\i}az, S. (2013). {\it Asymptotes and Perfect Curves}.
arxiv.org/abs/1307.6153.    Submitted to Computer Aided Geometric
Design.



 %\bibitem{fuster} Fuster, R.,  Gim\'enez, I.  (2008).
%{\it Variable Compleja y Ecuaciones Diferenciales}. Ed. Revert\'e.
\bibitem{conway} Conway, J.B.  (1995).
{\it Functions of One Complex Variable I}. Graduate Texts in Mathematics. Springer-Verlag. New York.


\bibitem{Duval89}  Duval, D.  (1989). {\it Rational Puiseux Expansion.}
Compositio Mathematica. Vol.    70, pp. 119--154.


%\bibitem{Fischer} Fischer, G.  (2001).
%{\it Plane Algebraic Curves}. Student Mathematical Library. Vol. 15. AMS.






\bibitem{Gao} Gao, B., Chen, Y.  (2012). {\it Finding the Topology of Implicitly defined two Algebraic Plane Curves.}
Journal of Systems Science and Complexity. Vol 25,  Issue 2, pp. 362-374.

\bibitem{lalo} Gonz\'alez-Vega, L.,  Necula, I. (2002). {\it Efficient Topology Determination of Implicitly defined Algebraic Plane Curves.} Comput. Aided Geom. Design. Vol. 19(9), pp. 719--743


\bibitem{Hong} Hong, H. (1996). {\it An Effective Method for Analyzing the Topology of Plane Real Algebraic Curves.}  Math. Comput. Simulation. Vol.
42,  pp. 572--582


\bibitem{HSW} Hoffmann, C.M., Sendra, J.R., Winkler, F. (1997).
{\it Parametric Algebraic Curves and Applications}. J. Symbolic
Computation. Vol. 23.

\bibitem{HL97} Hoschek, J., Lasser, D. (1993). {\it Fundamentals
of Computer Aided Geometric Design}.  A.K. Peters Wellesley MA.,
Ltd.



\bibitem{SWP}  Sendra, J.R., Winkler, F., P\'erez-D\'{\i}az, S. (2007). {\it Rational Algebraic Curves: A Computer Algebra Approach}. Series: Algorithms and Computation in Mathematics.  Vol. 22. Springer Verlag.



%\bibitem{spivak} Spivak,  M.  (2006). {\it Calculus.} Cambridge. Third edition. Vol. 1.

\bibitem{Stad00} Stadelmeyer,
P. (2000). {\it On the Computational Complexity of Resolving Curve
Singularities and Related Problems.} Ph.D. thesis, RISC-Linz, J.
Kepler Univ. Linz, Austria, Techn. Rep. RISC 00-31.


\bibitem{verger}  Verger-Gaugry, J-L. (2011). {\it Beta-Conjugates of Real Algebraic Numbers as
Puiseux Expansions}. Integers: Electronic Journal of Combinatorial
Number Theory. Proceedings of the Leiden Numeration Conference 2010.
Vol. 11B.

\bibitem{Walker} Walker, R.J.  (1950).  {\it Algebraic Curves}. Princeton University Press.



\bibitem{Zeng} Zeng, G. (2007). {\it Computing the Asymptotes for a Real Plane
Algebraic Curve.} Journal of Algebra. Vol. 316,  pp. 680–-705.

\end{thebibliography}
\end{document}